\documentclass{article}
\usepackage{amsfonts,latexsym,hyperref} 
\usepackage[ansinew]{inputenc}
\newcounter{satznum}
\newtheorem{theorem}{Theorem}[satznum]

\newtheorem{lemma}[theorem]{Lemma}

\newtheorem{proposition}[theorem]{Proposition}
\newenvironment{acknowledgement}
 {\begin{trivlist}\item[]{\bf Acknowledgement.}}
 {\end{trivlist}}
\newenvironment{remark}
 {\begin{trivlist}\item[]{\bf Remark.}}
 {\end{trivlist}}
\newenvironment{remarks}
 {\begin{trivlist}\item[]{\bf Remarks.}}
 {\end{trivlist}}
\newenvironment{example}
 {\begin{trivlist}\item[]{\bf Example.}}
 {\end{trivlist}}
\newenvironment{examples}
 {\begin{trivlist}\item[]{\bf Examples.}}
 {\end{trivlist}}
\newenvironment{proof}
 {\begin{trivlist}\item[]{\bf Proof.}}
 {\end{trivlist}}
{\makeatletter
\gdef\me{{\mathbb E}} 
\gdef\nz{{\mathbb N}} 
\gdef\pr{{\mathbb P}} 
\gdef\rz{{\mathbb R}} 
\gdef\gz{{\mathbb Z}} 
}
\makeatletter
\def\@MRExtract#1 #2!{#1}
\newcommand{\MR}[1]{
  \xdef\@MRSTRIP{\@MRExtract#1 !}
  \href{http://www.ams.org/mathscinet-getitem?mr=\@MRSTRIP}{MR\@MRSTRIP}}
\makeatother

\begin{document}
   \section*{ON THE BLOCK COUNTING PROCESS AND THE
   FIXATION LINE OF EXCHANGEABLE COALESCENTS}
   {\sc Florian Gaiser and Martin M\"ohle}\footnote{Mathematisches Institut, Eberhard Karls Universit\"at T\"ubingen,
   Auf der Morgenstelle 10, 72076 T\"ubingen, Germany, E-mail addresses: florian.gaiser@uni-tuebingen.de,
   martin.moehle@uni-tuebingen.de}
\begin{center}
   \today
\end{center}
\begin{abstract}
   We study the block counting process and the fixation line of
   exchangeable coalescents. Formulas for the infinitesimal rates of
   both processes are provided. It is shown that the block counting
   process is Siegmund dual to the fixation line. For exchangeable
   coalescents restricted to a sample of size $n$ and with dust we
   provide a convergence result for the block counting process as $n$
   tends to infinity. The associated limiting process is related to the
   frequencies of singletons of the coalescent. Via duality we obtain an
   analog convergence result for the fixation line of exchangeable
   coalescents with dust. The Dirichlet coalescent and the
   Poisson--Dirichlet coalescent are studied in detail.

   \vspace{2mm}

   \noindent Keywords: block counting process; coming down from infinity;
   Dirichlet coalescent; duality; dust;
   exchangeable coalescent; fixation line; Poisson--Dirichlet coalescent

   \vspace{2mm}

   \noindent 2010 Mathematics Subject Classification:
            Primary 60F05;   
            60J27            
            Secondary 92D15; 
            97K60            
\end{abstract}
\subsection{Introduction} \label{intro}
\setcounter{theorem}{0} Coalescent processes have attracted the
interest of many researchers over the last decades, mainly in
population genetics and probability. Most results in coalescent
theory concern coalescents with multiple collisions independently
introduced by Pitman \cite{pitman1} and Sagitov \cite{sagitov}.
Less is known for the full class of exchangeable coalescents
$\Pi=(\Pi_t)_{t\ge 0}$ allowing for simultaneous multiple
collisions of ancestral lineages. Note that $\Pi$ is a Markovian
process taking values in the space ${\cal P}$ of partitions of
$\nz:=\{1,2,\ldots\}$. Schweinsberg \cite{schweinsberg2} showed
that every exchangeable coalescent $\Pi$ can be characterized by a
finite measure $\Xi$ on the infinite simplex
\begin{equation} \label{delta}
   \Delta\ :=\
   \{x=(x_r)_{r\in\nz}\,:\,
   \mbox{$x_1\ge x_2\ge\cdots\ge 0,|x|:=\sum_{r=1}^\infty x_r\le 1$}\}.
\end{equation}
Exchangeable coalescents are therefore called $\Xi$-coalescents.
The aim of this article is to provide some more information on the
block counting process and the fixation line of the $\Xi$-coalescent and
on the relation between these two processes. We
therefore briefly recall the definition of the block counting
process and turn afterwards to the fixation line.

For $t\ge 0$ let $N_t$ denote the number of blocks of $\Pi_t$.
It is well known that $N:=(N_t)_{t\ge 0}$ is a Markovian process with state
space $S:=\nz\cup\{\infty\}$, called the block counting process of $\Pi$.
We use the notation $N_t^{(n)}$ for the number of blocks of
$\Pi_t$ restricted to a sample of size $n\in\nz$.
The block counting process has been studied intensively in the literature
with a main focus on coalescents with multiple collisions ($\Lambda$-coalescents).
We will revisit some of its properties throughout this article.

The definition of the fixation line is more involved. As in
\cite{schweinsberg2} decompose $\Xi=\Xi(\{0\})\delta_0+\Xi_0$ with
$\delta_0$ the Dirac measure at $0\in\Delta$ and
$\Xi_0$ having no atom at $0$. For $x=(x_r)_{r\in\nz}\in\Delta$
define $(x,x):=\sum_{r=1}^\infty x_r^2$ and $\nu({\rm
d}x):=\Xi_0({\rm d}x)/(x,x)$.
One possible definition of the fixation line is based on the lookdown
construction going back to Donnelly and Kurtz \cite{donnellykurtz1, donnellykurtz2}.
For some further information on the lookdown construction we refer the
reader to \cite{birknerblathmoehlesteinrueckentams}.
Imagine a population
consists of countably many individuals distinguished by their levels.
The level of an individual is a positive integer,
and the individual at time $t\ge 0$ at level $i\in\nz$ is denoted by $(t,i)$.
In the following we use Schweinsberg's Poisson process
construction \cite[Section 3]{schweinsberg2}.
Define $\nz_0:=\{0,1,2,\ldots\}$ and,
for every $x=(x_r)_{r\in\nz}\in\Delta$, let $P_x$ be the law of a sequence
$\xi=(\xi_1,\xi_2,\ldots)$ of independent and identically distributed
$\nz_0$-valued random variables with distribution
$\pr(\xi_1=0):=1-|x|$ and $\pr(\xi_1=r)=x_r$, $r\in\nz$.
Furthermore, for $i,j\in\nz$ with $i<j$
let $z_{ij}$ be the sequence $(z_1,z_2,\ldots)$ with
$z_i=z_j:=1$ and $z_k:=0$ for $k\notin\{i,j\}$.
Take a Poisson process on $[0,\infty)\times\nz_0^\nz$ with
intensity measure $\lambda\otimes \mu$, where $\lambda$ denotes Lebesgue
measure on $[0,\infty)$ and
\begin{equation} \label{mu}
   \mu(A)\ :=\ \Xi(\{0\})\sum_{{i,j\in\nz}\atop{i<j}}1_{\{z_{ij}\in A\}} + \int_\Delta P_x(A)\,\nu({\rm d}x)
\end{equation}
for all measurable $A\subseteq\nz_0^{\nz}$.
Each atom $(t,x)$ corresponds to a reproduction event which is defined
as follows. For $r\in\nz$ define $J_r:=\{j\in\nz:\xi_j=r\}$.
For $r\in\nz$ and $j\in J_r$, the individual $(t,j)$ is a child of
the individual $(t-,\min J_r)$. The other lineages are shifted upwards
keeping the order they had before the reproduction event. The construction
of this countable infinite population model is called the lookdown
construction.

We are now able to define the fixation line. Fix $i\in\nz$.
The levels of the offspring at time $t\ge 0$ of the individual $(0,i+1)$,
that is the individual at time $0$ at level $i+1$, form a subset of $\nz$,
whose minimal element (if it exists) we denote by $L_t^{(i)}+1$. If this
subset is empty, we define $L_t^{(i)}:=\infty$. For example, if
$J_1=\nz$ and $J_r=\emptyset$ for all $r\in\nz\setminus\{1\}$, then the
individual at time $0$ at level $2$ has no offspring at all, so in this
case we have $L_t^{(1)}=\infty$.

By construction, the process
$L^{(i)}:=(L_t^{(i)})_{t\ge 0}$ has state space $\{i,i+1,\ldots\}\cup\{\infty\}$
and non-decreasing paths. Moreover, $L_t^{(1)}\le L_t^{(2)}\le \cdots$.
We define $L_t:=L_t^{(1)}$ and $L:=(L_t)_{t\ge 0}$.
When $L_t$ reaches level $n$, all the individuals at time $t$ with
levels $1,\ldots,n$ are offspring of the single individual $(0,1)$,
an event called fixation in population genetics. The process $L$ is
hence called the fixation line. This process can be traced back to
Pfaffelhuber and Wakolbinger \cite{pfaffelhuberwakolbinger} for the
Kingman coalescent. For the $\Lambda$-coalescent the fixation line appears
in Labbé \cite{labbe} and was further studied by Hénard
\cite{henard1, henard2}.

We close the introduction by a brief summary of the organization of
the article. Section \ref{results} contains the
main results. Propositions \ref{ratesofn} and \ref{ratesofl} provide
formulas for the infinitesimal rates and the total rates of
the block counting process $N=(N_t)_{t\ge 0}$ and the fixation line
$L=(L_t)_{t\ge 0}$ for arbitrary $\Xi$-coalescents. Theorem \ref{duality}
shows that the block counting process $N$ is Siegmund dual to the fixation
line $L$. For $\Xi$-coalescents with dust,
Theorem \ref{theo1} provides a convergence result for $N$ and $L$ when
their initial state tends to infinity. The limiting processes are related
to the frequencies of singletons of the coalescent. In Sections
\ref{dirichlet} and \ref{poissondirichlet} we study the Dirichlet coalescent
and the Poisson--Dirichlet coalescent respectively. The proofs of the results
stated in Section \ref{results} are provided in Section \ref{proofs}.
The appendix deals with a duality relation for generalized Stirling numbers
being closely related to the Siegmund duality of $N$ and $L$.
\subsection{Results} \label{results}
\setcounter{theorem}{0}
Throughout the article we shall use the following subsets of the infinite
simplex $\Delta$ defined in (\ref{delta}). For $n\in\nz$ define
$\Delta_n:=\{x=(x_r)_{r\in\nz}\in\Delta\,:\,x_1+\cdots+x_n=1\}$.
Furthermore let $\Delta_f:=\bigcup_{n\in\nz}\Delta_n=\{x=(x_r)_{r\in\nz}\in\Delta\,:\,
x_1+\cdots+x_n=1\mbox{ for some }n\in\nz\}$ and
$\Delta^*:=\{x\in\Delta\,:\,|x|=1\}$. Note that
$\Delta_1\subset\Delta_2\subset\cdots\subset\Delta_f\subset\Delta^*\subset\Delta$
and that $\nu(\Delta_n)\le n\Xi(\Delta_n)<\infty$, since
$(x,x)=\sum_{r=1}^n x_r^2=1/n+\sum_{r=1}^n (x_r-1/n)^2\ge 1/n$ for
all $n\in\nz$ and $x\in\Delta_n$.

In order to state our first result it is convenient to introduce the
following urn model which is essentially a version of Kingman's paintbox
construction \cite[Section 8]{kingman}. Fix $x=(x_r)_{r\in\nz}\in\Delta$.
Recall that $|x|:=\sum_{i=1}^\infty x_i$ and define $x_0:=1-|x|$ for
convenience. Imagine a countable infinite number of boxes having labels
$r\in\nz_0$. Balls are allocated successively to these boxes where it is
assumed that every ball will go to box $r\in\nz_0$ with probability $x_r$
independently of the other balls. If $X_r(i,x)$ denotes the number of balls
in box $r\in\nz_0$ after $i\in\nz_0$ balls have been allocated, then clearly
$(X_0(i,x),X_1(i,x),X_2(i,x),\ldots)$ has an infinite multinomial
distribution with parameters $i$ and $(x_0,x_1,x_2,\ldots)$. The random
variable
\begin{equation} \label{y}
   Y(i,x)\ :=\ X_0(i,x) + \sum_{r=1}^\infty 1_{\{X_r(i,x)\ge 1\}}
\end{equation}
counts the balls contained in box $0$ plus the number of other
boxes which are non-empty. In the language of Kingman's paintbox,
when a ball goes to box $r\in\nz$, it will be painted with color $r$. The
box $0$ plays a distinguished role. Each ball going to box $0$ is
painted with a new color never seen before. $Y(i,x)$ is the number
of different colors after $i$ balls have been painted.

Proposition \ref{ratesofn} and \ref{ratesofl} below underline the well
known fact that the process $(Y(i,x))_{i\in\nz_0}$ plays a fundamental
role in coalescent theory. Proposition \ref{ratesofn} concerns the
infinitesimal rates of block counting process $N$. These rates are
essentially known from the literature (see, for example, \cite{freundmoehle}).
We state the result for the record and since the case $\Xi(\Delta_f)>0$
requires some attention.
\begin{proposition}[Rates of the block counting process] \label{ratesofn}
   Let $\Xi$ be a finite measure on $\Delta$ and let $\Pi$ be a $\Xi$-coalescent.
   The block counting process $N=(N_t)_{t\ge 0}$ of $\Pi$
   moves from state $i\in\nz$ to state $j\in\nz$ with $j<i$ at the rate
   \begin{equation} \label{qijgeneral}
      q_{ij}\ =\ \Xi(\{0\}) {i\choose 2}\delta_{j,i-1} +
      \int_\Delta \pr(Y(i,x)=j)\,\nu({\rm d}x)
   \end{equation}
   with $Y(i,x)$ defined via (\ref{y}).
   The probability below the integral in (\ref{qijgeneral})
   can be provided explicitly 
   as $\pr(Y(i,x)=j)=\sum_{k=1}^j f_{ijk}(x)$, where
   \begin{equation} \label{fijk}
      f_{ijk}(x)\ :=\ \frac{x_0^{j-k}}{(j-k)!}
      \sum_{{i_1,\ldots,i_k\in\nz}\atop{i_1+\cdots+i_k=i-j+k}}
      \frac{i!}{i_1!\cdots i_k!}
      \sum_{{r_1,\ldots,r_k\in\nz}\atop{r_1<\cdots<r_k}}
      x_{r_1}^{i_1}\cdots x_{r_k}^{i_k}.
   \end{equation}
   The total rates are 
   \begin{eqnarray}
      q_i
      & := & \sum_{j=1}^{i-1}q_{ij}
      \ = \ \Xi(\{0\}){i\choose 2} + \int_\Delta\pr(Y(i,x)<i)\,\nu({\rm d}x)\nonumber\\
      & = & \hspace{-1mm}\Xi(\{0\}){i\choose 2} + \int_\Delta
            \bigg(1-x_0^i-\sum_{k=1}^i {i\choose k}x_0^{i-k}
            \sum_{{r_1,\ldots,r_k\in\nz}\atop{\rm all\ distinct}} x_{r_1}\cdots x_{r_k}
            \bigg)
            \,\nu({\rm d}x),\ i\in\nz.\label{qigeneral}
   \end{eqnarray}
   Moreover, $q_{\infty j}=\nu(\Delta_j)-\nu(\Delta_{j-1})$
   for all $j\in\nz$ ($\Delta_0:=\emptyset$) and
   $q_{\infty\infty}=-\nu(\Delta_f)$.
\end{proposition}
\begin{remarks}
   \begin{enumerate}
      \item[1.] From $q_{i+1}-q_i=i\Xi(\{0\})+\int_\Delta\pr(Y(i+1,x)=i,Y(i,x)=i)\,\nu({\rm d}x)\ge 0$
         we conclude that $q_{i+1}\ge q_i$ with equality $q_{i+1}=q_i$ if
         and only if $\Xi(\Delta\setminus\Delta_{i-1})=0$, $i\in\nz$.
         Thus, if $\Xi(\Delta\setminus\Delta_f)>0$ then the
         total rates $q_i$, $i\in\nz$, are pairwise distinct.
      \item[2.] For the rates of the block counting process of the
         Dirichlet coalescent and the Poisson--Dirichlet coalescent
         we refer the reader to (\ref{diriqij1}), (\ref{diriqij2}) and
         (\ref{qijpd}).
      \item[3.] For the $\Lambda$-coalescent the rate (\ref{qijgeneral})
         reduces to the well known formula (see, for example, \cite{pitman1}
         or \cite[Eq.~(13)]{moehleewens})
         \begin{equation} \label{qij}
            q_{ij}
            \ =\ {i\choose{j-1}}\int_{[0,1]} x^{i-j-1}(1-x)^{j-1}\,\Lambda({\rm d}x),
            \qquad i,j\in\nz, j<i,
         \end{equation}
         and the total rates are \cite[Eq.~(14)]{moehleewens}
         \begin{eqnarray*}
            q_i
            & = & \Lambda(\{0\}){i\choose 2}
                  + \int_{(0,1]}\frac{1-(1-x)^i - ix(1-x)^{i-1}}{x^2}
                  \,\Lambda({\rm d}x),\qquad i\in\nz.
         \end{eqnarray*}
         For the $\beta(a,b)$-coalescent with parameters
         $a,b\in (0,\infty)$ the rate (\ref{qij}) reduces to
         \begin{equation} \label{qijbeta}
            q_{ij}\ =\ \frac{\Gamma(a+b)}{\Gamma(a)\Gamma(b)}
            \frac{\Gamma(i+1)}{\Gamma(i-2+a+b)}
            \frac{\Gamma(j-1+b)}{\Gamma(j)}
            \frac{\Gamma(i-j-1+a)}{\Gamma(i-j+2)},
            \qquad i,j\in\nz, j<i.
         \end{equation}
   \end{enumerate}
\end{remarks}
We now state the analog result for the rates of the fixation line
$L=(L_t)_{t\ge 0}$.
\begin{proposition}[Rates of the fixation line] \label{ratesofl}
   Let $\Xi$ be a finite measure on $\Delta$.
   The fixation line $L=(L_t)_{t\ge 0}$ moves from state $i\in\nz$
   to state $j\in\nz$ with $j>i$ at the rate
   \begin{equation} \label{gammaijgeneral}
      \gamma_{ij}\ =\ \Xi(\{0\}) {{j}\choose 2}\delta_{j,i+1} +
      \int_\Delta \pr(Y(j,x)=i,Y(j+1,x)=i+1)
      \,\nu({\rm d}x)
   \end{equation}
   with $Y(.,x)$ defined in (\ref{y}). Moreover, $\gamma_{i\infty}=\nu(\Delta_i)$
   for all $i\in\nz$ and
   $\gamma_{\infty\infty}=0$. The probability 
   below
   the integral in (\ref{gammaijgeneral}) can be provided explicitly,
   namely $\pr(Y(j,x)=i,Y(j+1,x)=i+1)=\sum_{k=1}^i g_{ijk}(x)$, where
   \begin{equation} \label{gijk}
      g_{ijk}(x)\ :=\
      \frac{x_0^{i-k}}{(i-k)!}
      \sum_{{i_1,\ldots,i_k\in\nz}\atop{i_1+\cdots+i_k=j-i+k}}
      \frac{j!}{i_1!\cdots i_k!}
      \sum_{{r_1,\ldots,r_k\in\nz}\atop{r_1<\cdots<r_k}}
      x_{r_1}^{i_1}\cdots x_{r_k}^{i_k}\bigg(1-\sum_{l=1}^k x_{r_l}\bigg).
   \end{equation}
   The total rates are $\gamma_i:=
   \sum_{j\in\{i+1,i+2,\ldots\}\cup\{\infty\}}\gamma_{ij}=q_{i+1}$, $i\in\nz$.
\end{proposition}
\begin{remarks}
   \begin{enumerate}
      \item[1.] In general $g_{ijk}(x)$ is not equal to $f_{jik}(x)$ (see
         (\ref{fijk})) because of the additional factor $1-\sum_{l=1}^k x_{r_l}$
         occurring on the right hand side in (\ref{gijk}). This
         additional factor comes from the fact that in the paintbox
         construction, on the event that box $0$ contains $i-k$ balls
         and that the boxes $r_1,\ldots,r_k$ are non-empty,
         $\{Y(j+1,x)=i+1\}$ corresponds to the event that ball $j+1$
         belongs to a box $r\in\nz_0\setminus\{r_1,\ldots,r_k\}$,
         which has probability $1-\sum_{l=1}^k x_{r_l}$.
      \item[2.] Note that $\gamma_{i+1}=\gamma_i$ if and only if
         $\Xi(\Delta\setminus\Delta_i)=0$, $i\in\nz$. If
         $\Xi(\Delta\setminus\Delta_f)>0$ then the total rates
         $\gamma_i$, $i\in\nz$, are pairwise distinct.
      \item[3.] For the $\Lambda$-coalescent the rate (\ref{gammaijgeneral})
         reduces to
         \begin{equation} \label{gammaij}
            \gamma_{ij}
            \ =\ {j\choose{i-1}}\int_{[0,1]}x^{j-i-1}(1-x)^i\,\Lambda({\rm d}x),
            \qquad i,j\in\nz, i<j,
         \end{equation}
         in agreement with \cite[Lemma 2.3]{henard2}. The total rates are
         $$
         \gamma_i\ =\ q_{i+1}
         \ =\ \Lambda(\{0\})\displaystyle{{i+1}\choose 2}
         +\int_{(0,1]} \frac{1-(1+ix)(1-x)^i}{x^2}\,\Lambda({\rm d}x),
         \qquad i\in\nz.
         $$
         For the $\Lambda$-coalescent the equality $\gamma_i=q_{i+1}$ was
         already observed by Hénard \cite{henard2}.

         For the Kingman coalescent ($\Lambda=\delta_0$) we have
         $\gamma_i=\gamma_{i,i+1}={{i+1}\choose 2}$ and
         $\gamma_{ij}=0$ for all $j\notin\{i,i+1\}$. Thus, $L$ is a pure
         birth process with state space $\nz\cup\{\infty\}$ and birth rates
         $\gamma_i={{i+1}\choose 2}$, $i\in\nz$. This process explodes,
         i.e. $\pr(L_t=\infty)>0$ for all $t>0$, since
         $\sum_{i=1}^\infty 1/\gamma_i<\infty$. In fact, $L$ reaches
         $\infty$ in finite time almost surely.
         For more details on explosion of $L$ we refer the reader to Remark
         3 after Theorem \ref{duality}.

         For the $\beta(a,b)$-coalescent with $a,b\in (0,\infty)$ the rates
         (\ref{gammaij}) of the fixation line process reduce to
         \begin{equation} \label{betarates}
            \gamma_{ij}
            \ =\
            \frac{\Gamma(a+b)}{\Gamma(a)\Gamma(b)}
            \frac{\Gamma(i+b)}{\Gamma(i)}
            \frac{\Gamma(j+1)}{\Gamma(j-1+a+b)}
            \frac{\Gamma(j-i-1+a)}{\Gamma(j-i+2)},\qquad i,j\in\nz, i<j.
         \end{equation}
         Simple formulas for the total rates $\gamma_i$ seem to be only
         available for particular parameter choices of $a$ and $b$. For
         example, for the $\beta(2-\alpha,\alpha)$-coalescent with
         parameter $0<\alpha<2$, which has attracted the interest of
         several researchers, the fixation line has total rates
         \begin{eqnarray*}
            \gamma_i
            & = & \sum_{j=i+1}^\infty\gamma_{ij}
            \ = \ \frac{1}{\Gamma(2-\alpha)\Gamma(\alpha)}
                  \frac{\Gamma(i+\alpha)}{\Gamma(i)}
                  \sum_{j=i+1}^\infty\frac{\Gamma(j-i+1-\alpha)}{\Gamma(j-i+2)}\\
            & = & \frac{1}{\Gamma(2-\alpha)\Gamma(\alpha)}
                  \frac{\Gamma(i+\alpha)}{\Gamma(i)}\frac{\Gamma(2-\alpha)}{\alpha}
            \ = \ \frac{\Gamma(i+\alpha)}{\Gamma(\alpha+1)\Gamma(i)}
            \ = \ \prod_{k=1}^{i-1} \frac{k+\alpha}{k},\quad i\in\nz.
         \end{eqnarray*}
         In particular, $\gamma_i=i$ for the Bolthausen--Sznitman coalescent
         ($\alpha=1$). Another class of beta-coalescents for which a nice
         formula for the total rate $\gamma_i$ is available is the
         $\beta(3,b)$-coalescent with parameter $b>0$. In this case the fixation
         line process has total rates
         \begin{eqnarray*}
            \gamma_i
            & = & \sum_{j=i+1}^\infty \gamma_{ij}
            \ = \ \frac{\Gamma(3+b)}{\Gamma(3)\Gamma(b)}
                  \frac{\Gamma(i+b)}{\Gamma(i)}
                  \sum_{j=i+1}^\infty \frac{\Gamma(j+1)}{\Gamma(j+b+2)}\\
            & = & \frac{\Gamma(3+b)}{\Gamma(3)\Gamma(b)}
                  \frac{\Gamma(i+b)}{\Gamma(i)}
                  \frac{\Gamma(i+2)}{b\Gamma(i+b+2)}\\
            & = & \frac{(b+1)(b+2)}{2}\frac{i(i+1)}{(i+b)(i+b+1)},
                  \qquad i\in\nz.
         \end{eqnarray*}
   \end{enumerate}
\end{remarks}
We now turn to the duality of the block counting process $N$ and the
fixation line $L$. The following result (Theorem \ref{duality}) is a
reformulation and generalization of Lemma 2.1 of Hénard \cite{henard2}.
Note that Theorem \ref{duality} holds for any $\Xi$-coalescent.
Recall the notation $S:=\nz\cup\{\infty\}$.
\begin{theorem}[Siegmund duality of N and L] \label{duality}
   Let $\Pi$ be a $\Xi$-coalescent and let $N=(N_t)_{t\ge 0}$ and
   $L=(L_t)_{t\ge 0}$ denote
   the block counting process and the fixation line of $\Pi$ respectively.
   Then $N$ is dual (in the sense of
   Liggett \cite[p.~84, Definition 3.1]{liggett}) to
   $L$ with respect to the Siegmund duality kernel
   $H:S^2\to\{0,1\}$ defined via $H(i,j):=1$ for $i\le j$ and $H(i,j):=0$
   otherwise, i.e.
   \begin{eqnarray*}
      \pr(N_t^{(i)}\le j)
      & = & \pr(N_t\le j\,|\,N_0=i)
      \ = \ \me(H(N_t,j)\,|\,N_0=i)\\
      & = & \me(H(i,L_t)\,|\,L_0=j)
      \ = \ \pr(L_t\ge i\,|\,L_0=j)
      \ =\ \pr(L_t^{(j)}\ge i)
   \end{eqnarray*}
   for all $i,j\in S$ and all $t\ge 0$. If $Q=(q_{ij})_{i,j\in S}$ and
   $\Gamma=(\gamma_{ij})_{i,j\in S}$ denote the generator matrices of
   $N$ and $L$ respectively then $q_{i,\le j}:=\sum_{k\in S,k\le j}q_{ik}
   =\sum_{k\in S,k\ge i}\gamma_{jk}=:\gamma_{j,\ge i}$ for all $i,j\in S$.
\end{theorem}
\begin{remarks}
   \begin{enumerate}
      \item[1.] Note that $L$ is dual to $N$ with respect to the transposed
         duality kernel $H^\top:S^2\to\{0,1\}$ defined via $H^\top(i,j):=
         H(j,i)$ for all $i,j\in S$.
         Let $\tilde{S}:=S\setminus\{1,\infty\}=\nz\setminus\{1\}$ denote
         the set of states which are transient for $N$ and $L$. Define the
         matrix $\tilde{H}:=(H(i,j))_{i,j\in\tilde{S}}$.
         The matrix inverse $\tilde{H}^{-1}$ of $\tilde{H}$ has entries
         $$
         \tilde{H}^{-1}(i,j)
         \ =\ \left\{
            \begin{array}{cl}
               1 & \mbox{if $j=i$,}\\
               -1 & \mbox{if $j=i+1$,}\\
               0 & \mbox{otherwise.}
            \end{array}
         \right.
         $$
         For $i,j\in\tilde{S}$ define $g_{ij}:=\int_0^\infty \pr(N_t^{(i)}=j)
         \,{\rm d}t$ and $\widehat{g}_{ij}:=\int_0^\infty\pr(L_t^{(i)}=j)
         \,{\rm d}t$, and let $G=(g_{ij})_{i,j\in\tilde{S}}$ and
         $\widehat{G}=(\widehat{g}_{ij})_{i,j\in\tilde{S}}$ denote the Green matrices
         of $N$ and $L$ respectively (see, for example, Norris
         \cite[p.~145]{norris}).

From Theorem \ref{duality} it follows that
         $G\tilde{H}=\tilde{H}\widehat{G}^\top$.
         Thus, $G$ can be
         computed from $\widehat{G}$ and vice versa via
         $G=\tilde{H}\widehat{G}^\top\tilde{H}^{-1}$ and $\widehat{G}
         =(\tilde{H}^{-1}G\tilde{H})^\top$, i.e.
         $$
         g_{i2}\ =\ \sum_{{k\in\tilde{S}}\atop{k\ge i}}\widehat{g}_{2k},
         \qquad
         g_{ij}
         \ =\ \sum_{{k\in\tilde{S}}\atop{k\ge i}}(\widehat{g}_{jk}-\widehat{g}_{j-1,k}),
         \qquad i\in\tilde{S}, j\in\tilde{S}\setminus\{2\}
         $$
         and
         $$
         \widehat{g}_{ij}\ =\ \sum_{{k\in\tilde{S}}\atop{k\le i}}(g_{jk}-g_{j+1,k}),
         \qquad i,j\in\tilde{S}.
         $$
      \item[2.] As in \cite{schweinsberg2} we say that $\Pi$ comes down
         from infinity if $\pr(N_t<\infty)=1$ for all $t>0$. We say that
         $\Pi$ stays infinite if $\pr(N_t=\infty)=1$ for all $t>0$. Note
         that there exist coalescents, for example the Dirichlet coalescent
         studied in Section \ref{dirichlet}, that neither come down from
         infinity nor stay infinite. We refer the reader to
         \cite{schweinsberg1, schweinsberg2} and \cite{herrigermoehle} for
         methods to determine wether a coalescent $\Pi$ comes
         down from infinity or stays infinite.
      \item[3.] We say that $L$ does not explode if $T_\infty:=
         \inf\{t>0\,:\,L_{t-}=\infty\}=\infty$ almost surely. Note that
         $L$ does not explode if and only if $\pr(L_t<\infty)=1$ for all
         $t\ge 0$. By the general explosion criterium for Markov chains,
         $L$ does not explode if and only if $\sum_{n=0}^\infty
         1/\gamma_{\chi_n}=\infty$ almost surely, where
         $\chi=(\chi_n)_{n\in\nz_0}$ denotes the jump chain
         of $L$. Note that $\chi$ has transition probabilities
         $p_{ij}:=\gamma_{ij}/\gamma_i$, $1\le i<j\le\infty$.
      \item[4.] If $\Pi$ comes down from infinity then $L$ explodes.
         Moreover, $\Pi$ stays infinite if and only if $L$ does not explode.

         Proof. If $\Pi$ comes down from infinity, then $\Pi$ eventually
         becomes absorbed almost surely, i.e. $\lim_{t\to\infty}\pr(N_t=1)=1$.
         Thus, by duality, $\pr(L_t=\infty)=\pr(N_t=1)>0$ for all sufficiently large $t$,
         i.e. $L$ explodes.

         If $\Pi$ stays infinite, then $\pr(N_t=\infty)=1$ for all
         $t\ge 0$. Thus, $\pr(L_t=\infty)=\pr(N_t=1)=0$ for all $t\ge 0$, i.e.
         $L$ does not explode. Conversely, suppose that
         $L$ does not explode. Then, by the first statement, $\Pi$ does not
         come down from infinity. Moreover, we must have $\gamma_{i\infty}=0$
         for all $i\in\nz$, because otherwise every $L_t$ would be equal to $\infty$
         with positive probability. Since $\gamma_{i\infty}=\nu(\Delta_i)$ it
         follows that $\nu(\Delta_f)=\lim_{i\to\infty}\nu(\Delta_i)=0$.
         Thus, $\Xi(\Delta_f)=0$. But under the additional assumption that
         $\Xi(\Delta_f)=0$ the coalescent does not come down
         from infinity if and only if the coalescent stays infinite
         (see \cite{schweinsberg2}).
      \item[5.] Recall that $L$ explodes if $\Pi$ comes down
         from infinity. The converse holds under the additional assumption
         that $\Xi(\Delta_f)=0$, but it does not hold in general. Examples where
         $L$ explodes but $\Pi$ does not come down from infinity are provided
         in Section \ref{dirichlet} (Dirichlet coalescent).
   \end{enumerate}
\end{remarks}
A $\Xi$-coalescent $\Pi=(\Pi_t)_{t\ge 0}$ has proper frequencies if, for
all times $t\ge 0$, the frequency of singletons $S_t$ of $\Pi_t$ satisfies
$S_t=0$ almost surely. For a precise definition of $S_t$ we refer the reader to
\cite[Section 3]{moehleproper}. Schweinsberg \cite[Proposition 30]{schweinsberg2}
showed that $\Pi$ does not have proper frequencies
if and only if $\Xi(\{0\})=0$ and $\int_\Delta |x|\,\nu({\rm d}x)<\infty$.
In this case the process $Z:=(Z_t)_{t\ge 0}:=(-\log S_t)_{t\ge 0}$ (with the
convention $-\log 0:=\infty$) is a drift-free
subordinator with state space $[0,\infty]$ and Laplace exponent
\begin{equation} \label{phi}
   \Phi(\eta)\ =\ \int_\Delta
   (1-(1-|x|)^\eta)\,\nu({\rm d}x),\qquad \eta\ge 0.
\end{equation}
Note that $\me(S_t^\eta)=\me(e^{-\eta Z_t})=e^{-t\Phi(\eta)}$, $\eta\ge 0$.
A coalescent without proper frequencies is also called a coalescent with
dust component \cite{gnediniksanovmarynych}. Theorem \ref{theo1} below clarifies
the asymptotic behavior of the block counting process $N^{(n)}$ and the
fixation line $L^{(n)}$ as $n\to\infty$ for $\Xi$-coalescents
with dust. Note that we use the conventions $e^{-\infty}:=0$, $e^\infty:=\infty$
and $1/0:=\infty$.
\begin{theorem} \label{theo1}
   Let $\Pi$ be a $\Xi$-coalescent with dust, i.e. $\Xi(\{0\})=0$
   and $\int_\Delta |x|\,\nu({\rm d}x)<\infty$. Then, the following two
   assertions hold.
   \begin{enumerate}
      \item[a)] As $n\to\infty$ the scaled block counting process
         $(N_t^{(n)}/n)_{t\ge 0}$ converges in $D_{[0,1]}[0,\infty)$
         to the frequency of singleton process $S=(S_t)_{t\ge 0}=(e^{-Z_t})_{t\ge 0}$.
      \item[b)] As $n\to\infty$ the scaled fixation line process
         $(L_t^{(n)}/n)_{t\ge 0}$ converges in $D_{[1,\infty]}[0,\infty)$
         to the inverse frequency of singleton process $(1/S_t)_{t\ge 0}=
         (e^{Z_t})_{t\ge 0}$.
   \end{enumerate}
\end{theorem}
\begin{remarks}
   \begin{enumerate}
      \item[1.] Theorem \ref{theo1} can be stated logarithmically as follows.
         For $\Xi$-coalescents with dust, as $n\to\infty$, both processes
         $(\log n-\log N_t^{(n)})_{t\ge 0}$ and $(\log L_t^{(n)}-\log n)_{t\ge 0}$
         converge in $D_{[0,\infty]}[0,\infty)$ to the drift-free subordinator
         $Z$ with Laplace exponent (\ref{phi}).
         Clearly, Theorem \ref{theo1} holds for
         $\Lambda$-coalescents with dust, for example for
         $\beta(a,b)$-coalescents with $a>1$ and $b>0$.
      \item[2.] If $\Xi$ is concentrated on $\Delta^*$ then the
         coalescent has dust if and only if $\nu$ is finite.
         In this case the Laplace exponent (\ref{phi}) satisfies
         $\Phi(\eta)=\nu(\Delta^*)$, $\eta>0$, so $S_t$
         has the same distribution as $1_{\{T_f>t\}}$ for all $t\ge 0$,
         where $T_f$ is exponentially distributed with parameter
         $\nu(\Delta^*)$. Examples are the Dirichlet coalescent
         studied in Section \ref{dirichlet}, the
         Poisson--Dirichlet coalescent \cite[Section 3]{sagitov2} and
         the two-parameter Poisson--Dirichlet coalescent
         \cite[Section 6]{moehleproper}. More information on the
         two-parameter Poisson--Dirichlet coalescent is provided
         in Section \ref{poissondirichlet}.
      \item[3.] Theorem \ref{theo1} excludes dust-free coalescents.
         For the Bolthausen--Sznitman coalescent we refer the
         reader to \cite{kuklamoehle}, \cite{moehlemittag} and
         \cite{schweinsberg3} for asymptotic results concerning the
         block counting process $N^{(n)}$ and the fixation line $L^{(n)}$.
         For dust-free coalescents different from the Bolthausen--Sznitman
         coalescent we leave the asymptotic analysis of $N^{(n)}$ and
         $L^{(n)}$ for future work.
   \end{enumerate}
\end{remarks}
\subsection{The Dirichlet coalescent} \label{dirichlet}
\setcounter{theorem}{0}
   Let $X:=(X_1,\ldots,X_N)$ be symmetric Dirichlet distributed with
   parameters $N\in\nz$ and $\alpha>0$ and let
   $X_{(1)}\ge \cdots\ge X_{(N)}$ denote the order statistics of
   $X$.
   We consider the $\Xi$-coalescent when the characteristic measure $\nu$
   is the distribution of
   $(X_{(1)},\ldots,X_{(N)},0,0,\ldots)$. We call this coalescent
   the Dirichlet coalescent with parameters $N\in\nz$ and $\alpha>0$.
   Note that $\nu$ is concentrated on $\Delta_N$.
   The Dirichlet coalescent neither comes down from
   infinity nor stays infinite. More precisely, using the notation $A_f$
   and $T_f$ from \cite{schweinsberg2}, $\pr(N_t=\infty)=\pr(T_f>t)=
   e^{-t}$ for all $t\ge 0$, since (see \cite[Lemma 41]{schweinsberg2})
   $T_f$ is exponentially distributed with parameter $\int_\Delta
   P_x(A_f)\,\nu({\rm d}x)=\int_{\Delta_f} P_x(A_f)\,\nu({\rm d}x)=
   \nu(\Delta_f)=1$. Note that the fixation line $L=(L_t)_{t\ge 0}$
   explodes, since the coalescent does not stay infinite.

   In agreement with \cite{hsushiue} we use the notation
   $[x|y]_n:=\prod_{k=0}^{n-1}(x+ky)$ and $(x|y)_n:=\prod_{k=0}^{n-1}(x-ky)$
   for $x,y\in\rz$ and $n\in\nz_0$ with the convention that empty products are
   equal to $1$. We furthermore write $[x]_n:=[x|1]_n$ and $(x)_n:=(x|1)_n$.
   The proof of the following lemma is given at the end of this section.
\begin{lemma}[Rates of the block counting process] \label{diriratesofn}
   The block counting process of the
   Dirichlet coalescent with parameters $N\in\nz$ and $\alpha>0$ has
   infinitesimal rates
   \begin{equation} \label{diriqij1}
      q_{ij}
      \ =\ \frac{(N\alpha|\alpha)_j}{[N\alpha]_i}\,S_\alpha(i,j)
      \qquad 1\le j<i,
   \end{equation}
   where the $S_\alpha(i,j):=S(i,j;-1,\alpha,0)$ are the generalized
   Stirling numbers defined in \cite{hsushiue} satisfying the recursion
   $S_\alpha(i+1,j)=S_\alpha(i,j-1)+(i+\alpha j)S_\alpha(i,j)$. Alternatively,
   \begin{equation} \label{diriqij2}
      q_{ij}\ =\ \frac{\displaystyle{N\choose j}}{\displaystyle{{N\alpha+i-1}\choose i}}
      \sum_{{i_1,\ldots,i_j\in\nz}\atop{i_1+\cdots+i_j=i}}
      {{i_1+\alpha-1}\choose{i_1}}\cdots{{i_j+\alpha-1}\choose{i_j}},
      \qquad 1\le j<i.
   \end{equation}
\end{lemma}
\begin{remarks}
   \begin{enumerate}
      \item[1.] The Dirichlet coalescent is closely related to the Chinese
         restaurant process \cite{pitman2}. Imagine a restaurant with
         $N\in\nz$ tables each of infinite capacity. Customers enter
         successively the restaurant. When the $(i+1)$th customer arrives and
         $j$ tables are already occupied (by at least one person), the
         customer sits at an empty table with probability
         $(N\alpha-j\alpha)/(N\alpha+i)$. This corresponds to Pitman's
         \cite{pitman2} Chinese restaurant process with $\kappa:=\alpha$
         and $\theta:=N\alpha$. Let $K_i$ denote the number of occupied
         tables after the $i$th customer has been seated. It is easily
         verified by induction on $i$ that $K_i$ has distribution
         $\pr(K_i=j)=(N\alpha|\alpha)_jS_\alpha(i,j)/[N\alpha]_i$, $j\in\{1,\ldots,i\}$. Note that $K_i$ has mean
         $$
         \me(K_i)\ =\ N-N\frac{[(N-1)\alpha]_i}{[N\alpha]_i}.
         $$
      \item[2.] The Dirichlet coalescent has total rates
         $q_i:=\sum_{j=1}^{i-1}q_{ij}
         =\sum_{j=1}^{i-1}\pr(K_i=j)=1-\pr(K_i=i)
         =1-(N\alpha|\alpha)_i/[N\alpha]_i$.
         Note that $0=q_1<q_2<\cdots<q_N<1=q_{N+1}=q_{N+2}=\cdots$. The
         Dirichlet coalescent therefore serves as an example
         showing that in general the total rates of a coalescent do not need
         to be pairwise distinct.
   \end{enumerate}
\end{remarks}
\begin{examples}
   \begin{enumerate}
      \item[1.] For $\alpha=1$ the Stirling number $S_1(i,j)$ is the Lah
         number $S(i,j;-1,1,0)=\frac{i!}{j!}{{i-1}\choose{j-1}}$ and we conclude
         that $q_{ij}={N\choose j}{{i-1}\choose{j-1}}/{{N+i-1}\choose i}$.
         In this case $K_i$ is hypergeometric distributed with parameters
         $N+i-1$, $N$ and $i$. The total rates are
         $q_i=1-N!(N-1)!/(N-i)!/(N+i-1)!$ for $i\le N$ and $q_i=1$ for
         $i>N$.
      \item[2.] For $\alpha\to\infty$ it follows that
         $q_{ij}=N^{-i}(N)_jS(i,j)$, where the $S(i,j)$ are the usual Stirling
         numbers of the second kind. In this case $K_i$ counts the number of
         non-empty boxes when $i$ balls are allocated at random to $N$ boxes.
         This corresponds to the Dirac $\Xi$-coalescent where the measure
         $\nu$ assigns its total mass $1$ to the single point
         $x\in\Delta$ whose first $N$ coordinates are all equal to $1/N$.
      \item[3.] For $\alpha\to 0$ and $N\to\infty$ such that
         $N\alpha\to\theta\in (0,\infty)$ the rates $q_{ij}$ converge to
         those of the Poisson--Dirichlet coalescent with parameter
         $\theta$ and $\alpha=0$ studied in the following Section
         \ref{poissondirichlet}.
   \end{enumerate}
\end{examples}
In the following we provide the asymptotics of some functionals of
the Dirichlet $n$-coalescent when the sample size $n$ tends to infinity.
By Theorem \ref{theo1} and Remark 3 thereafter, $(N_t^{(n)}/n)_{t\ge 0}$
converges in $D_{[0,1]}[0,\infty)$ to $(S_t)_{t\ge 0}$ as
$n\to\infty$ and $(L_t^{(n)}/n)_{t\ge 0}$ converges in
$D_{[1,\infty]}[0,\infty)$ to $(1/S_t)_{t\ge 0}$, where
$S_t:=1_{\{T_f>t\}}$ and $T_f$ is exponentially distributed with
parameter $1$.

Let $C_n$ denote the number of jumps and
$\tau_n:=\inf\{t>0\,:\,N_t=1\}$ denote the absorption time
of the Dirichlet $n$-coalescent.
The following lemma clarifies the asymptotics of $C_n$ and $\tau_n$ as
$n\to\infty$. Its proof is given at the end of this section.
\begin{lemma}[Asymptotics of the number of jumps and the absorption time] \label{asy}
\ \\
   For the Dirichlet coalescent with parameter $N\in\nz$ and
   $\alpha>0$ the following two statements hold.
   \begin{enumerate}
      \item[(i)] The number of jumps $C_n$ converges to $C_\infty:=1+C_N$ in distribution as $n\to\infty$.
         The limit $C_\infty$ has distribution
         $\pr(C_\infty=k)=r_{N1}^{(k-1)}$, $1\le k\le N$, and
         $r_{ij}:=q_{ij}/q_i$, $1\le j<i$, denote the transition probabilities of the
         jump chain of the block counting process.
      \item[(ii)] The absorption time $\tau_n$ converges to $\tau_\infty:=E+\tau_N$ in distribution as
         $n\to\infty$, where $E$ is standard exponentially distributed
         and independent of $\tau_N$.
   \end{enumerate}
\end{lemma}
We now turn to the fixation line $L=(L_t)_{t\ge 0}$ of the Dirichlet
coalescent.
\begin{lemma}[Rates of the fixation line] \label{diriratesofl}
   The fixation line of the Dirichlet coalescent with
   parameters $N\in\nz$ and $\alpha>0$ has rates
   \begin{equation} \label{dirigammaij}
      \gamma_{ij}
      \ =\ \frac{(N\alpha|\alpha)_{i+1}}{[N\alpha]_{j+1}}S_\alpha(j,i),
      \qquad i,j\in\nz, i<j.
   \end{equation}
   Moreover $\gamma_{i\infty}=0$ for $i\in\{1,\ldots,N-1\}$ and
   $\gamma_{i\infty}=1$ for $i\in\{N,N+1,\ldots\}$.
\end{lemma}
\begin{remark}
   Note that $\gamma_{ij}=\frac{N\alpha-i\alpha}{N\alpha+j}q_{ji}=
   \pr(K_j=i,K_{j+1}=i+1)=\pr(K_j\le i,K_{j+1}>i)=\pr(K_j\le i)-\pr(K_{j+1}\le i)$.
   Summation over all $j\in\{i+1,i+2,\ldots\}\cup\{\infty\}$ shows that
   the fixation line has total rates $\gamma_i=q_{i+1}$, $i\in\nz$, in
   agreement with Proposition \ref{ratesofl}.
\end{remark}
\begin{example}
   For $\alpha=1$ we obtain
   $$
   \gamma_{ij}\ =\ \frac{N-i}{N+j}q_{ji}
   \ =\ {{N-1}\choose i}{{j-1}\choose{j-i}}\bigg/{{N+j}\choose N}.
   $$
\end{example}
By duality (Theorem \ref{duality}), for all $i,j\in\nz$ with $i>j$, the
two quantities
\begin{equation} \label{diri1}
   q_{i,\le j}\ =\ \sum_{k=1}^j q_{ik}
   \ =\ \frac{1}{[N\alpha]_i}\sum_{k=1}^j (N\alpha|\alpha)_j S_\alpha(i,j)
\end{equation}
and
\begin{equation} \label{diri2}
   \gamma_{j,\ge i}
   \ =\ \sum_{{k\in\nz\cup\{\infty\}}\atop{k\ge i}}\gamma_{jk}
   \ =\
   \left\{
      \begin{array}{cl}
         (N\alpha|\alpha)_{j+1}\displaystyle\sum_{k=i}^\infty \frac{S_\alpha(k,j)}{[N\alpha]_{k+1}} & \mbox{if $j<N$}\\
         \gamma_{j\infty}=1 & \mbox{if $j\ge N$}
      \end{array}
   \right.
\end{equation}
coincide. The equality of (\ref{diri1}) and (\ref{diri2}) follows alternatively
from \cite[Lemma 4.1]{moehlegreenwood}, applied to the Markov chain $K$, or
from Lemma \ref{applemma} in the appendix, applied with $a:=-1$, $b:=\alpha$,
$r:=0$ and $t:=N\alpha$. Note that
$\lim_{i\to\infty}q_{i,\le j}=\lim_{i\to\infty}\pr(K_i\le j)=0$ for all
$j<N$, since all states $j<N$ of the Markov chain $K$ are transient.

\vspace{2mm}

In the remaining part of this section we prove
Lemma \ref{diriratesofn}, Lemma \ref{asy} and Lemma \ref{diriratesofl}.

\begin{proof} (of Lemma \ref{diriratesofn})
   Since the measure $\nu$ is concentrated on $\Delta_N$ it follows from
   (\ref{qijgeneral}) and (\ref{fijk}) that
   \begin{equation} \label{local}
      q_{ij}
      \ = \ \frac{i!}{j!}
            \sum_{{i_1,\ldots,i_j\in\nz}\atop{i_1+\cdots+i_j=i}}
            \frac{\phi(i_1,\ldots,i_j)}{i_1!\cdots i_j!},
   \end{equation}
   where
   $$
   \phi(i_1,\ldots,i_j)
   \ :=\ \int_\Delta\sum_{{r_1,\ldots,r_j\in\nz}\atop{\rm all\;distinct}}
   x_{r_1}^{i_1}\cdots x_{r_j}^{i_j}\,\nu({\rm d}x)
   \ =\ \int_\Delta \sum_{{r_1,\ldots,r_j=1}\atop{\rm all\;distinct}}^N
   x_{r_1}^{i_1}\cdots x_{r_j}^{i_j}\,\nu({\rm d}x).
   $$
   The function below the latter integral is symmetric with respect to
   $x_1,\ldots,x_N$. Thus,
   \begin{eqnarray*}
      \phi(i_1,\ldots,i_j)
      & = & \int \sum_{{r_1,\ldots,r_j=1}\atop{\rm all\;distinct}}^N
            x_{r_1}^{i_1}\cdots x_{r_j}^{i_j}
            \,D_N(\alpha)({\rm d}x_1,\ldots,{\rm d}x_N)
      \ = \ (N)_j\me(X_1^{i_1}\cdots X_j^{i_j}),
   \end{eqnarray*}
   where the last equality holds, since the Dirichlet distribution
   $D_N(\alpha)$ is symmetric and, hence, the integrals over
   each summand is the same.
   The moments of the symmetric Dirichlet distribution $D_N(\alpha)$ are well
   known (see, for example, \cite[p.~488]{kotzbalakrishnanjohnson}) to be
   $\me(X_1^{i_1}\cdots X_j^{i_j})=[\alpha]_{i_1}\cdots [\alpha]_{i_j}/[N\alpha]_i$,
   where $i:=i_1+\cdots+i_j$.
   Plugging all this into (\ref{local}) leads to
   \begin{eqnarray*}
      q_{ij}
      & = & \frac{i!}{j!}\frac{(N)_j}{[N\alpha]_i}
            \sum_{{i_1,\ldots,i_j\in\nz}\atop{i_1+\cdots+i_j=i}}
            \frac{[\alpha]_{i_1}\cdots[\alpha]_{i_j}}{i_1!\cdots i_j!}\\
      & = & \frac{{N\choose j}}{{{N\alpha+i-1}\choose i}}
            \sum_{{i_1,\ldots,i_j\in\nz}\atop{i_1+\cdots+i_j=i}}
            {{i_1+\alpha-1}\choose{i_1}}\cdots{{{i_j+\alpha-1}\choose{i_j}}},
   \end{eqnarray*}
   which is (\ref{diriqij2}). Since $(N\alpha|\alpha)_j=
   (N)_j\alpha^j$ and
   $$
   S_\alpha(i,j)\ =\ \frac{i!}{j!}
   \sum_{{i_1,\ldots,i_j\in\nz}\atop{i_1+\cdots+i_j=i}}
   \frac{1}{i_1!\cdots i_j!}\frac{\Gamma(\alpha+i_1)}{\Gamma(\alpha+1)}
   \cdots\frac{\Gamma(\alpha+i_j)}{\Gamma(\alpha+1)}
   $$
   it follows that the right hand side of (\ref{diriqij2}) is equal to
   the right hand side of (\ref{diriqij1}).\hfill$\Box$
\end{proof}
\begin{proof} (of Lemma \ref{asy})
   Clearly, $(C_n)_ {n\in\nz}$ satisfies the distributional recursion
   $C_n\stackrel{d}{=}1+C_{I_n}$ with initial condition $C_1=0$, where
   $I_n$ denotes the state of the jump chain of the block counting process
   of the Dirichlet $n$-coalescent with parameters $N\in\nz$ and $\alpha>0$
   after its first jump. By Lemma \ref{diriratesofn}
   and the remarks thereafter, $I_n$ has distribution
   $$
   \pr(I_n=k)
   \ =\ \frac{q_{nk}}{q_n}
   \ =\ \frac{\pr(K_n=k)}{1-\pr(K_n=n)},\qquad 1\le k<n.
   $$
   Since $K_n\to N$ in distribution as $n\to\infty$ we conclude that
   $I_n\to N$ in distribution as $n\to\infty$. Thus,
   $C_n\stackrel{d}{=}1+C_{I_n}\to 1+C_N$ in distribution as $n\to\infty$.

   The sequence $(\tau_n)_{n\in\nz}$ satisfies the distributional recursion
   $\tau_n\stackrel{d}{=}E_n+\tau_{I_n}$ with initial condition $\tau_1=0$,
   where $E_n$ is independent of $I_n$ and exponentially distributed with parameter
   $q_n=1-\pr(K_n=n)$. Since $q_n\to 1$
   as $n\to\infty$ and $I_n\to N$ in distribution as $n\to\infty$ we
   conclude that $\tau_n\stackrel{d}{=}E_n+\tau_{I_n}\to E+\tau_N$
   as $n\to\infty$, where $E$ is standard exponentially distributed
   and independent of $\tau_N$.\hfill$\Box$
\end{proof}
\begin{proof} (of Lemma \ref{diriratesofl})
   Since $\nu$ is concentrated on $\Delta_N$ it follows from
   (\ref{gammaijgeneral}) and (\ref{gijk}) that
   $$
   \gamma_{ij}\ =\
   \frac{j!}{i!}\sum_{{j_1,\ldots,j_i\in\nz}\atop{j_1+\cdots+j_i=j}}
   \frac{\psi(j_1,\ldots,j_i)}{j_1!\cdots j_i!},
   $$
   where
   $$
   \psi(j_1,\ldots,j_i)
   \ :=\ \int_\Delta
   \sum_{{r_1,\ldots,r_i=1}\atop{\rm all\;distinct}}^N
   x_{r_1}^{j_1}\cdots x_{r_i}^{j_i}\bigg(1-\sum_{k=1}^i x_{r_k}\bigg)
   \,\nu({\rm d}x).
   $$
   The same arguments as in the proof of Lemma \ref{diriratesofn}
   show that
   \begin{eqnarray*}
      \psi(j_1,\ldots,j_i)
      & = & (N)_i \me\bigg(X_1^{j_1}\cdots X_i^{j_i}\bigg(1-\sum_{k=1}^i X_k\bigg)\bigg)\\
      & = & (N)_i\bigg(
               \frac{[\alpha]_{j_1}\cdots [\alpha]_{j_i}}{[N\alpha]_j}
               - \sum_{k=1}^i
               \frac{[\alpha]_{j_1}\cdots [\alpha]_{j_k+1}\cdots[\alpha]_{j_i}}{[N\alpha]_{j+1}}
            \bigg)\\
      & = & (N)_i \frac{[\alpha]_{j_1}\cdots[\alpha]_{j_i}}{[N\alpha]_{j+1}}
            \bigg((N\alpha+j)-\sum_{k=1}^i(\alpha+j_k)\bigg)\\
      & = & \frac{N\alpha-i\alpha}{N\alpha+j}(N)_i\frac{[\alpha]_{j_1}\cdots[\alpha]_{j_i}}{[N\alpha]_j}
      \ = \ \frac{N\alpha-i\alpha}{N\alpha+j}\phi(j_1,\ldots,j_i).
   \end{eqnarray*}
   It follows that $\gamma_{ij}=q_{ji}(N\alpha-i\alpha)/(N\alpha+j)$ and
   (\ref{dirigammaij}) follows from (\ref{diriqij1}). Moreover, by
   Proposition \ref{ratesofl}, $\gamma_{i\infty}=\nu(\Delta_i)$ for all
   $i\in\nz$. It remains to note that $\nu(\Delta_i)=0$ for $i<N$ and
   $\nu(\Delta_i)=1$ for $i\ge N$.\hfill$\Box$
\end{proof}
\subsection{The Poisson--Dirichlet coalescent} \label{poissondirichlet}
\setcounter{theorem}{0}
   If $\Xi(\{0\})=0$ and if the measure $\nu({\rm d}x)=\Xi_0({\rm d}x)/(x,x)$ is the
   Poisson--Dirichlet distribution with parameters $0\le\alpha<1$
   and $\theta>-\alpha$ then the $\Xi$-coalescent is called the
   two-parameter Poisson--Dirichlet coalescent \cite{moehleproper}.
   Note that $\nu$ is concentrated on $\Delta^*\setminus\Delta_f$. Since
   $\int_\Delta |x|\,\nu({\rm d}x)=\nu(\Delta^*)=1<\infty$ it follows that
   this coalescent has dust and hence cannot come down from infinity.
   Since $\Xi(\Delta_f)=0$, this coalescent stays infinite, and, hence,
   $L$ does not explode. The associated block counting process has rates
   (see \cite[p.~2170]{moehleproper})
   \begin{equation} \label{qijpd}
      q_{ij}\ =\ c_{j,\alpha,\theta}\frac{\Gamma(\theta+\alpha j)}{\Gamma(\theta+i)}s_\alpha(i,j),
      \qquad i>j,
   \end{equation}
   where
   $$
   c_{j,\alpha,\theta}\ :=\ \prod_{k=1}^j \frac{\Gamma(\theta+1+(k-1)\alpha)}{\Gamma(1-\alpha)\Gamma(\theta+k\alpha)}
   $$
   and
   $$
   s_\alpha(i,j)\ :=\ \frac{i!}{j!}\sum_{{i_1,\ldots,i_j\in\nz}\atop{i_1+\cdots+i_j=i}}
   \frac{\Gamma(i_1-\alpha)\cdots\Gamma(i_j-\alpha)}{i_1!\cdots i_j!}
   $$
   is a kind of generalized absolute Stirling number of the first kind
   satisfying the recursion $s_\alpha(i+1,j)=\Gamma(1-\alpha)s_\alpha(i,j-1)+
   (i-\alpha j)s_\alpha(i,j)$.
   More precisely, $s_\alpha(i,j)/(\Gamma(1-\alpha))^j$ coincides
   with the generalized Stirling number $S(i,j;-1,-\alpha,0)$ as defined in
   Hsu and Shiue \cite{hsushiue}.
   For $\alpha=0$ (and hence $\theta>0$) the rate $q_{ij}$ reduces to
   $$
   q_{ij}\ =\ \theta^j\frac{\Gamma(\theta)}{\Gamma(\theta+i)} s(i,j),
   \qquad i>j,
   $$
   where the $s(i,j):=S(i,j;-1,0,0)$ are the (usual) absolute Stirling numbers of the first kind.
   For $\theta=0$ (and hence $0<\alpha<1$) we obtain
   $$
   q_{ij}\ =\ \frac{\alpha^{j-1}}{(\Gamma(1-\alpha))^j}\frac{(j-1)!}{(i-1)!}\,s_\alpha(i,j)
   \ =\ \alpha^{j-1}\frac{(j-1)!}{(i-1)!}S(i,j;-1,-\alpha,0),\qquad i>j.
   $$
   In order to compute the total rates of the block counting process
   of the two-parameter Poisson--Dirichlet coalescent we proceed as follows.
   Let $K=(K_n)_{n\in\nz_0}$ be a
   Markov chain with state space $\nz_0$, $K_0:=0$, $K_1:=1$ and
   transition probabilities $p_k(n):=\pr(K_{n+1}=k+1\,|\,K_n=k)
   :=(\theta+\alpha k)/(\theta+n)$ and $\pr(K_{n+1}=k\,|\,K_n=k)=1-p_k(n)$
   for $n\in\nz$ and $k\in\{1,\ldots,n\}$. Note that
   $1\le K_n\le n$, $n\in\nz$.
   As for the Dirichlet coalescent one may interpret $K_n$ as the number of
   occupied tables in a particular Chinese restaurant process. When
   customer $n+1$ enters the restaurant and $k$ tables are already occupied,
   he sits at an empty table with probability $p_k(n)$. In the following it
   is verified by induction on $n\in\nz$ that $K_n$ has distribution
   (see also Pitman \cite[p.~65, Eq.~(3.11)]{pitman2})
   $$
   \pr(K_n=k)\ =\ c_{k,\alpha,\theta}
   \,\frac{\Gamma(\theta+\alpha k)}{\Gamma(\theta+n)}
   \,s_\alpha(n,k),\qquad k\in\nz_0.
   $$
   For $n=1$ this is obvious, since $K_1=1$, $c_{1,\alpha,\theta}=\Gamma(\theta+1)/
   \Gamma(1-\alpha)/\Gamma(\theta+\alpha)$ and $s_\alpha(1,1)=\Gamma(1-\alpha)$.
   The induction step from $n\in\nz$ to $n+1$ works as follows. By the
   Markov property,
   $\pr(K_{n+1}=k)=p_{k-1}(n)\,\pr(K_n=k-1) + (1-p_k(n))\,\pr(K_n=k)$.
   By induction,
   \begin{eqnarray*}
      p_{k-1}(n)\,\pr(K_n=k-1)
      & = & \frac{\theta+\alpha(k-1)}{\theta+n}
            \,c_{k-1,\alpha,\theta}
            \,\frac{\Gamma(\theta+\alpha(k-1))}{\Gamma(\theta+n)}
            \,s_\alpha(n,k-1)\\
      & = & \frac{\Gamma(\theta+1+\alpha(k-1))}{\Gamma(\theta+n+1)}
            \,c_{k-1,\alpha,\theta}
            \,s_\alpha(n,k-1)\\
      & = & \frac{\Gamma(\theta+k\alpha)\Gamma(1-\alpha)}{\Gamma(\theta+n+1)}
            \,c_{k,\alpha,\theta}\,s_\alpha(n,k-1)
   \end{eqnarray*}
   and
   \begin{eqnarray*}
      (1-p_k(n))\,\pr(K_n=k)
      & = & \frac{n-\alpha k}{\theta+n}\,
            c_{k,\alpha,\theta}\,\frac{\Gamma(\theta+k\alpha)}{\Gamma(\theta+n)}
            s_\alpha(n,k)\\
      & = & \frac{\Gamma(\theta+k\alpha)}{\Gamma(\theta+n+1)}
            \,c_{k,\alpha,\theta}\,(n-\alpha k)\,s_\alpha(n,k).
   \end{eqnarray*}
   Summation of these two terms yields
   \begin{eqnarray*}
      \pr(K_{n+1}=k)
      & = & \frac{\Gamma(\theta+\alpha k)}{\Gamma(\theta+n+1)}\,c_{k,\alpha,\theta}
            \,\big(\Gamma(1-\alpha)s_\alpha(n,k-1)+(n-\alpha k)s_\alpha(n,k)\big)\\
      & = & \frac{\Gamma(\theta+\alpha k)}{\Gamma(\theta+n+1)}
            \,c_{k,\alpha,\theta}\,s_\alpha(n+1,k),
   \end{eqnarray*}
   which completes the induction.
   As a consequence, the block counting process of the two-parameter
   Poisson--Dirichlet coalescent has total rates
   \begin{eqnarray*}
      q_i
      & = & \sum_{j=1}^{i-1} q_{ij}
      \ = \ \sum_{j=1}^{i-1}\pr(K_i=j)
      \ =\ 1-\pr(K_i=i)\\
      & = & 1-c_{i,\alpha,\theta}\frac{\Gamma(\theta+\alpha i)}{\Gamma(\theta+i)}s_\alpha(i,i)\\
      & = & 1-\frac{\Gamma(\theta+\alpha i)}{\Gamma(\theta+i)}
            \prod_{k=1}^i \frac{\Gamma(\theta+1+(k-1)\alpha)}{\Gamma(\theta+k\alpha)} ,
   \qquad i\in\nz.
   \end{eqnarray*}
   For $\alpha=0$ the total rates reduce to
   $q_i=1-\theta^i\Gamma(\theta)/\Gamma(\theta+i)$, $i\in\nz$, and for
   $\theta=0$ we obtain $q_i=1-\alpha^{i-1}$, $i\in\nz$.

   Let us now turn to the fixation line. The rates $\gamma_{ij}$, $i<j$,
   of the fixation line are obtained similarly as the rates $q_{ij}$ as
   follows. Since the measure $\nu$ is concentrated on $\Delta^*$
   it follows from (\ref{gammaijgeneral}) and (\ref{gijk}) that
   $$
   \gamma_{ij}
   \ =\ \frac{j!}{i!}\sum_{{j_1,\ldots,j_i\in\nz}\atop{j_1+\cdots+j_i=j}}
   \frac{I(j_1,\ldots,j_i)}{j_1!\cdots j_i!},
   $$
   where
   $$
   I(j_1,\ldots,j_i)\ :=\ \int_\Delta
   \sum_{{r_1,\ldots,r_i\in\nz}\atop{\rm all\ distinct}}
   x_{r_1}^{j_1}\cdots x_{r_i}^{j_i}\bigg(1-\sum_{k=1}^ix_{r_k}\bigg)
   \nu({\rm d}x).
   $$
   By \cite[Eq.~(2.1)]{handa},
   $I(j_1,\ldots,j_i)
   =\int_{\rz^i} x_1^{j_1}\cdots x_i^{j_i}(1-\sum_{k=1}^ix_k)
   \,\mu_i({\rm d}x_1,\ldots,{\rm d}x_i),
   $
   where $\mu_i$ denotes the $i$th correlation measure associated with the
   Poisson--Dirichlet distribution. The density (correlation function)
   of $\mu_i$ is known explicitly (see, for example, \cite[Theorem 2.1]{handa})
   and we obtain
   $
   I(j_1,\ldots,j_i)
   \ =\ c_{i,\alpha,\theta}\int_{\Delta_i} x_1^{j_1-\alpha-1}\cdots x_i^{j_i-\alpha-1}(1-\sum_{k=1}^i x_k)^{\theta+\alpha i}
   \,{\rm d}x_1\cdots{\rm d}x_i.
   $
   The last integral is known (Liouville's integration formula), and it
   follows that
   $$
   I(j_1,\ldots,j_i)\ =\ c_{i,\alpha,\theta}
   \frac{\Gamma(j_1-\alpha)\cdots\Gamma(j_i-\alpha)\Gamma(\theta+\alpha i+1)}{\Gamma(\theta+j+1)},
   $$
   $j:=j_1+\cdots+j_i>i$. Plugging this expression into the above formula for
   $\gamma_{ij}$ leads to
   \begin{eqnarray}
      \gamma_{ij}
      & = & c_{i,\alpha,\theta}
            \frac{\Gamma(\theta+\alpha i+1)}{\Gamma(\theta+j+1)}
           \frac{j!}{i!}
           \sum_{{j_1,\ldots,j_i\in\nz}\atop{j_1+\cdots+j_i=j}}
           \frac{\Gamma(j_1-\alpha)\cdots\Gamma(j_i-\alpha)}{j_1!\cdots j_i!}\nonumber\\
      & = & c_{i,\alpha,\theta}
            \frac{\Gamma(\theta+\alpha i+1)}{\Gamma(\theta+j+1)}s_\alpha(j,i),
            \qquad i<j. \label{gammaijpd}
   \end{eqnarray}
   For $\alpha=0$ the rate $\gamma_{ij}$ reduces to
   $$
   \gamma_{ij}\ = \ \theta^i\frac{\Gamma(\theta+1)}{\Gamma(\theta+j+1)}s(j,i),
   \qquad i<j
   $$
   whereas for $\theta=0$ we obtain
   $$
   \gamma_{ij}\ =\ \frac{\alpha^i}{(\Gamma(1-\alpha))^i}
   \frac{i!}{j!}s_\alpha(j,i)
   \ =\ \alpha^i \frac{i!}{j!}S(j,i;-1,-\alpha,0),\qquad i<j.
   $$
   In particular,
   \begin{equation} \label{pd1}
   q_{i,\le j}
   \ =\ \sum_{k=1}^j q_{ik}
   \ =\ \frac{1}{\Gamma(\theta+i)}\sum_{k=1}^j c_{k,\alpha,\theta}\Gamma(\theta+\alpha k) s_\alpha(i,k),\qquad i>j,
   \end{equation}
   and
   \begin{equation} \label{pd2}
   \gamma_{j,\ge i}
   \ =\ \sum_{k=i}^\infty \gamma_{jk}
   \ =\ c_{j,\alpha,\theta}\Gamma(\theta+\alpha j+1)
        \sum_{k=i}^\infty \frac{s_\alpha(k,j)}{\Gamma(\theta+k+1)},
   \qquad i>j.
   \end{equation}
   The two expressions (\ref{pd1}) and (\ref{pd2}) are equal by duality
   (Theorem \ref{duality}). The equality of (\ref{pd1}) and
   (\ref{pd2}) also follows from \cite[Lemma 4.1]{moehlegreenwood},
   applied to the Markov chain $K$. Alternatively, one may apply Lemma
   \ref{applemma} in the appendix with $a:=-1$, $b:=-\alpha$, $r:=0$ and $t:=\theta$.
   Note however that formally the case $\theta=0$ is not covered by
   Lemma \ref{applemma}. As a last option one may prove the
   equality of (\ref{pd1}) and (\ref{pd2}) directly (using the recursion
   for generalized Stirling numbers) as follows. We have
   \begin{eqnarray*}
      &   & \hspace{-1cm}q_{ki}-q_{k+1,i}\\
      & = &
            c_{i,\alpha,\theta}
            \,\frac{\Gamma(\theta+\alpha i)}{\Gamma(\theta+k)}\,s_\alpha(k,i)
            -
            c_{i,\alpha,\theta}\,\frac{\Gamma(\theta+\alpha i)}{\Gamma(\theta+k+1)}\,s_\alpha(k+1,i)\\
      & = & c_{i,\alpha,\theta}\frac{\Gamma(\theta+\alpha i)}{\Gamma(\theta+k)}s_\alpha(k,i)\\
      &   & \hspace{2cm}
            - c_{i,\alpha,\theta}
            \frac{\Gamma(\theta+\alpha i)}{\Gamma(\theta+k+1)}
            \big(\Gamma(1-\alpha)s_\alpha(k,i-1)+(k-\alpha i)s_\alpha(k,i)\big)\\
      & = & c_{i,\alpha,\theta}\frac{\Gamma(\theta+\alpha i)}{\Gamma(\theta+k+1)}
            s_\alpha(k,i)((\theta+k)-(k-\alpha i))\\
      &   & \hspace{2cm}
            - c_{i,\alpha,\theta}\frac{\Gamma(\theta+\alpha i)}{\Gamma(\theta+k+1)}
            \Gamma(1-\alpha) s_\alpha(k,i-1)\\
      & = & c_{i,\alpha,\theta}\,\frac{\Gamma(\theta+1+\alpha i)}{\Gamma(\theta+k+1)}\,s_\alpha(k,i)
            -
            c_{i-1,\alpha,\theta}\frac{\Gamma(\theta+1+(i-1)\alpha)}{\Gamma(\theta+k+1)}\,
            s_\alpha(k,i-1).
   \end{eqnarray*}
   Summation over all $i\in\{1,\ldots,j\}$ yields
   $$
   q_{k,\le j} - q_{k+1,\le j}
   \ =\ c_{j,\alpha,\theta}\frac{\Gamma(\theta+1+\alpha j)}{\Gamma(\theta+k+1)}s_\alpha(k,j).
   $$
   Another summation over all $k\ge i$ yields
   $q_{i,\le j}=c_{j,\alpha,\theta}\Gamma(\theta+1+\alpha j)
   \sum_{k=i}^\infty s_\alpha(k,j)/\Gamma(\theta+k+1)$, which shows that
   (\ref{pd1}) and (\ref{pd2}) coincide.
\subsection{Proofs} \label{proofs}
\setcounter{theorem}{0}
\begin{proof} (of Proposition \ref{ratesofn})
   The formulas (\ref{qijgeneral}) and (\ref{qigeneral}) for the rates
   $q_{ij}$ and the total rates $q_i$ of the block counting process are
   known from the literature \cite[Eqs.~(1.2) and (1.3)]{freundmoehle}.
   For given $x=(x_r)_{r\in\nz}\in\Delta$
   the block counting process jumps from $\infty$ to $j\in\nz$ if and only if
   $x_1+\cdots+x_j=1$ and $x_1,\ldots,x_j>0$, i.e. if and only if
   $x\in\Delta_j\setminus\Delta_{j-1}$ with the convention
   $\Delta_0:=\emptyset$. Integration with respect to $\nu$
   yields $q_{\infty j}=\nu(\Delta_j\setminus\Delta_{j-1})=\nu(\Delta_j)-
   \nu(\Delta_{j-1})$ for all $j\in\nz$.
   Since the generator $Q=(q_{ij})_{i,j\in S}$ is conservative, it follows
   that $q_{\infty\infty}=-\sum_{j\in\nz}q_{\infty j}
   =-\sum_{j\in\nz}(\Delta(\nu_j)-\Delta(\nu_{j-1}))
   =-\lim_{n\to\infty}\nu(\Delta_n)
   =-\nu(\Delta_f)$. 
   \hfill$\Box$
\end{proof}
\begin{proof} (of Proposition \ref{ratesofl})
   We generalize Hénard's proof of Lemma 2.3 in \cite{henard2}.
   Recall the pathwise definition of the fixation line based on the lookdown
   construction provided in the introduction.
   Assume first that $\Xi(\{0\})=0$.
   The fixation line jumps from $i\in\nz$ to $j\in\nz$ with $j>i$ if
   and only if there exists
   $k\in\{1,\ldots,i\}$ and $1\le r_1<\cdots<r_k$ such that
   \begin{enumerate}
        \item[(i)] exactly $i-k$ of the individuals $1,\ldots,j$ belong
         to $J_0:=\bigcup_{r\in\nz}J_r$,
      \item[(ii)] for every $l\in\{1,\ldots,k\}$ at least one of the
         individuals $1,\ldots,j$ belongs to $J_{r_l}$ and
      \item[(iii)] the individual $j+1$ does not belong to
         $J_{r_1}\cup\cdots\cup J_{r_k}$.
   \end{enumerate}
   For fixed $x=(x_r)_{r\in\nz}\in\Delta$ this event has probability
   $$
   {j\choose{i-k}}(1-|x|)^{i-k}
   \sum_{{i_1,\ldots,i_k\in\nz}\atop{i_1+\cdots+i_k=j-(i-k)}}
   \frac{(j-(i-k))!}{i_1!\cdots i_k!}
   x_{r_1}^{i_1}\cdots x_{r_k}^{i_k}
   \bigg(1-\sum_{l=1}^k x_{r_l}\bigg).
   $$
   Summing this probability over all $k\in\{1,\ldots,i\}$ and
   $1\le r_1<\cdots<r_k$ and integrating with respect to the law $\nu$
   yields $\gamma_{ij}=\int_\Delta\sum_{k=1}^i g_{ijk}(x)\,\nu({\rm d}x)$
   with $g_{ijk}(x)$ defined in (\ref{gijk}).
   From the definition of $Y(.,x)$ in (\ref{y}) it follows that $\sum_{k=1}^i g_{ijk}(x)
   =\pr(Y(j,x)=i,Y(j+1,x)=i+1)$.

   If $\Xi(\{0\})>0$ then the rate
   $\gamma_{ij}$ increases by $\Xi(\{0\}){j\choose 2}\delta_{j,i+1}$, since
   ${j\choose 2}\delta_{j,i+1}$ is the rate at which the fixation line of the
   Kingman coalescent jumps from $i$ to $j$. Thus (\ref{gammaijgeneral}) is
   established. Similarly, given $x=(x_r)_{r\in\nz}\in\Delta$,
   the fixation line jumps from $i\in\nz$ to $\infty$ if and only if
   $x_1+\cdots+x_i=1$, i.e. if and only if $x\in\Delta_i$. Integration
   with respect to $\nu$ yields $\gamma_{i\infty}=\nu(\Delta_i)$, $i\in\nz$.
   Clearly, $\gamma_{\infty\infty}=0$, since the state $\infty$ is absorbing.

   It remains to determine the total rates $\gamma_i$, $i\in\nz$. From
   the definition of $Y(.,x)$ in (\ref{y}) via the
   paintbox construction it follows that $Y(j+1,x)-Y(j,x)\in\{0,1\}$ for
   all $j\in\nz$ and $x\in\Delta$. Thus, for all $i,j\in\nz$,
   $\pr(Y(j,x)=i,Y(j+1,x)=i+1)=\pr(Y(j,x)\le i,Y(j+1,x)>i)
   =\pr(Y(j,x)\le i)-\pr(Y(j+1,x)\le i)$. Summation over
   all $j\in\nz$ with $j>i$ yields
   \begin{eqnarray*}
      \sum_{{j\in\nz}\atop{j>i}}\pr(Y(j,x)=i,Y(j+1,x)=i+1)
      & = & \sum_{{j\in\nz}\atop{j>i}} \big(\pr(Y(j,x)\le i) - \pr(Y(j+1,x)\le i)\big)\\
      & = & \pr(Y(i+1,x)\le i) - \lim_{k\to\infty}\pr(Y(k,x)\le i)\\
      & = & \left\{
            \begin{array}{cl}
               \pr(Y(i+1,x)\le i) & \mbox{if $x\in\Delta\setminus\Delta_i$,}\\
               0 & \mbox{if $x\in\Delta_i$,}
            \end{array}
            \right.
   \end{eqnarray*}
   since $\pr(Y(k,x)\le i)\to 0$ as $k\to\infty$ if $x\in\Delta\setminus\Delta_i$
   and $\pr(Y(k,x)\le i)=1$ for all $k\in\nz$ if $x\in\Delta_i$.
   Thus, the total rates of the fixation line are
   \begin{eqnarray*}
      \gamma_i
      & = & \gamma_{i\infty} + \sum_{{j\in\nz}\atop{j>i}} \gamma_{ij}\\
      & = & \gamma_{i\infty} + \Xi(\{0\}) {{i+1}\choose 2} + \int_{\Delta\setminus\Delta_i}\pr(Y(i+1,x)\le i)\,\nu({\rm d}x)\\
      & = & \Xi(\{0\}) {{i+1}\choose 2} + \int_\Delta \pr(Y(i+1,x)\le i)\,\nu({\rm d}x),
      \qquad i\in\nz.
   \end{eqnarray*}
   A comparison with the total rate $q_i$ of the block counting process
   shows that $\gamma_i=q_{i+1}$, $i\in\nz$. \hfill$\Box$
\end{proof}
\begin{proof} (of Theorem \ref{duality})
   From $Y(k+1,x)-Y(k,x)\in\{0,1\}$ for all $k\in\nz$ and $x\in\Delta$
   we conclude that
   \begin{eqnarray*}
      \pr(Y(k,x)=j,Y(k+1,x)=j+1)
      & = & \pr(Y(k,x)\le j,Y(k+1,x)>j)\\
      & = & \pr(Y(k,x)\le j) - \pr(Y(k+1,x)\le j)
   \end{eqnarray*}
   for all $j,k\in\nz$ and $x\in\Delta$.
   Integration with respect to $\nu$ and taking the formula
   (\ref{gammaijgeneral}) for the rates of the fixation line into account,
   it follows for all $j,k\in\nz$ with $j<k$ that
   \begin{eqnarray*}
      \gamma_{jk}
      & = & \Xi(\{0\}) {k\choose 2}\delta_{k,j+1} +
            \int_\Delta \pr(Y(k,x)=j,Y(k+1,x)=j+1)\,\nu({\rm d}x)\\
      & = & \Xi(\{0\}) {k\choose 2}\delta_{j,k-1} +
            \int_\Delta
            \big(\pr(Y(k,x)\le j) - \pr(Y(k+1,x)\le j)\big)\,\nu({\rm d}x)\\
      & = & \sum_{l=1}^j \bigg(
            \Xi(\{0\}) {k\choose 2}\delta_{l,k-1}
            + \int_\Delta\pr(Y(k,x)=l)\,\nu({\rm d}x)
            \\
      &   & \hspace{2cm}
            - \Xi(\{0\}) {{k+1}\choose 2}\delta_{lk}
            - \int_\Delta\pr(Y(k+1,x)=l)\,\nu({\rm d}x)
            \bigg)\\
      & = & \sum_{l=1}^j (q_{kl} - q_{k+1,l})
      \ = \ q_{k,\le j} - q_{k+1,\le j},\qquad j,k\in\nz,j<k,
   \end{eqnarray*}
   where the second last equality holds by (\ref{qijgeneral}).
   Let $i,j\in\nz$ with $i>j$. Summing over all $k\in\nz$ with $k\ge i$ yields
   $$
   \sum_{{k\in\nz}\atop{k\ge i}}\gamma_{jk}
   \ =\ \sum_{{k\in\nz}\atop{k\ge i}}(q_{k,\le j}-q_{k+1,\le j})
   \ =\ q_{i,\le j} - \lim_{k\to\infty} q_{k,\le j}
   \ =\ q_{i,\le j} - \nu(\Delta_j).
   $$
   The last equality holds since, by (\ref{qijgeneral}),
   \begin{eqnarray*}
      q_{k,\le j}
      & = & \sum_{l=1}^j q_{kl}
      \ = \ \Xi(\{0\}) {k\choose 2}\delta_{j,k-1}+ \int_\Delta \pr(Y(k,x)\le j)\,\nu({\rm d}x)\\
      & \to & \int_\Delta 1_{\Delta_j}(x)\,\nu({\rm d}x)
      \ = \ \nu(\Delta_j)
   \end{eqnarray*}
   as $k\to\infty$ by dominated convergence. Note that
   $\pr(Y(k,x)\le j)\le\pr(Y(j+1,x)\le j)$ for all $k>j$ and that
   the dominating map $x\mapsto\pr(Y(j+1,x)\le j)$ is $\nu$-integrable.

   Since $\gamma_{j\infty}=\nu(\Delta_j)$ it follows that
   \begin{equation} \label{gendual}
      q_{i,\le j}\ =\ \sum_{{k\in S}\atop{k\ge i}}\gamma_{jk}
      \ =\ \gamma_{j,\ge i}
   \end{equation}
   for all $i,j\in\nz$ with $i>j$. Eq.~(\ref{gendual}) holds as well for
   $i,j\in\nz$ with $i\le j$ since in this case both sides
   in (\ref{gendual}) are equal to $0$. Moreover,
   $q_{i,\le\infty}=0=\gamma_{\infty,\ge i}$ for all $i\in S$ and
   $q_{\infty,\le j}=\nu(\Delta_j)=\gamma_{j\infty}=\gamma_{j,\ge\infty}$
   for all $j\in\nz$.
   Thus, (\ref{gendual}) holds for all $i,j\in S$.

   Let $Q=(q_{ij})_{i,j\in S}$ and
   $\Gamma=(\gamma_{ij})_{i,j\in S}$ denote the generator matrices of $N$
   and $L$ respectively and let $H=(h_{ij})_{i,j\in S}$ denote the matrix
   with entries $h_{ij}:=1$ for $i\le j$ and $h_{ij}:=0$ for $i>j$. Since
   $(QH)_{ij}
   =\sum_{k\in S, k\le j} q_{ik}h_{kj}
   =\sum_{k\in S, k\le j} q_{ik}
   =q_{i,\le j}$
   and
   $(H\Gamma^\top)_{ij}
   =\sum_{k\in S, k\ge i} h_{ik} \gamma_{jk}
   =\sum_{k\in S, k\ge i}\gamma_{jk}
   =\gamma_{j,\ge i}
   $
   it follows that $QH=H\Gamma^\top$.
   It follows that $Q^kH=H(\Gamma^\top)^k$ for all $k\in\nz$ and, hence,
   $e^{tQ}H=H(e^{t\Gamma})^\top$ for all $t\ge 0$. Since $e^{tQ}$ and
   $e^{t\Gamma}$ are the transition matrices of the block counting process
   $N=(N_t)_{t\ge 0}$ and the fixation line $L=(L_t)_{t\ge 0}$ respectively,
   this shows that $N$ is Siegmund
   dual to $L$ with respect
   to the kernel $H$.\hfill$\Box$
\end{proof}
\begin{remark}
   For $\Lambda$-coalescents
   Lemma 2.1 of Hénard \cite{henard2} essentially states the
   Siegmund duality of $N$ and $L$ 
   and Lemma 2.4 of
   \cite{henard2} is a reformulation of this duality in terms of
   the generators of $N$ and $L$.
\end{remark}
\begin{proof} (of Theorem \ref{theo1})
   For $n\in\nz$ and $i\in\{1,\ldots,n\}$ let $\tau_{n,i}:=\inf\{t>0\,:\,
   \mbox{$i$ is not a singleton of $\Pi_t^{(n)}$}\}$ denote the length of
   the $i$th external branch of $\Pi^{(n)}$. For every $i\in\nz$ the
   sequence $(\tau_{n,i})_{n\ge i}$ is non-increasing in $n$ with
   $\tau_{n,i}\searrow\tau_i$ almost surely as $n\to\infty$, where
   $\tau_i:=\inf\{t>0\,:\,\mbox{$\{i\}$ is not a singleton of $\Pi_t$}\}$
   denotes the length of the $i$th external branch of $\Pi$. Let
   $t_1,\ldots,t_k\ge 0$. Conditional on $S_{t_1},\ldots,S_{t_k}$,
   the probability that $i$ is still a singleton at time $t_i$ for all
   $i\in\{1,\ldots,k\}$, is $S_{t_1}\cdots S_{t_k}$. Thus,
   $\pr(\tau_1>t_1,\ldots,\tau_k>t_k)=\me(S_{t_1}\cdots S_{t_k})$.
   In particular, $\pr(\tau_1>t,\ldots,\tau_k>t)=\me(S_t^k)$.

   \vspace{2mm}

   {\bf Proof of part a).}
   For $n\in\nz$ and $t\ge 0$
   decompose $N_t^{(n)}=E_t^{(n)}+I_t^{(n)}$, where $E_t^{(n)}:=
   \sum_{i=1}^n 1_{\{\tau_{n,i}>t\}}$ and $I_t^{(n)}:=N_t^{(n)}-E_t^{(n)}$
   denotes the number of singleton and non-singleton blocks of $\Pi_t^{(n)}$
   respectively. We think of $E_t^{(n)}$ and $I_t^{(n)}$ as the number of
   `external' and `internal' blocks of $\Pi_t^{(n)}$ and proceed similar
   as in the proof of Theorem 3 of \cite{moehleproper}. For
   $t\ge 0$ and $n,k\in\nz$,
   \begin{eqnarray*}
      \me((E_t^{(n)})^k)
      & = & \me((1_{\{\tau_{n,1}>t\}}+\cdots+1_{\{\tau_{n,n}>t\}})^k)\\
      & = & \sum_{{k_1,\ldots,k_n\in\nz_0}\atop{k_1+\cdots+k_n=k}}
            \frac{k!}{k_1!\cdots k_n!}
            \me(1_{\{\tau_{n,1}>t\}}^{k_1}\cdots 1_{\{\tau_{n,n}>t\}}^{k_n})\\
      & = & \sum_{j=1}^k {n\choose j}
            \sum_{{k_1,\ldots,k_j\in\nz}\atop{k_1+\cdots+k_j=k}}
            \me(1_{\{\tau_{n,1}>t\}}^{k_1}\cdots 1_{\{\tau_{n,j}>t\}}^{k_j}),
   \end{eqnarray*}
   where the last equality holds
   since the random variables $\tau_{n,i}$, $i\in\{1,\ldots,n\}$, are
   exchangeable. Thus,
   \begin{equation} \label{moments}
      \me((E_t^{(n)})^k)
      \ = \ \sum_{j=1}^k (n)_j S(k,j)\pr(\tau_{n,1}>t,\ldots,\tau_{n,j}>t),
            \qquad t\ge 0,n,k\in\nz,
   \end{equation}
   where $(n)_j:=n(n-1)\cdots(n-j+1)$ and $S(.,.)$ denote the Stirling
   numbers of the second kind. Dividing by $n^k$, letting $n\to\infty$
   and noting that $(n)_j/n^k\to\delta_{jk}$ it follows that
   \begin{equation} \label{momconv}
      \lim_{n\to\infty}\me((E_t^{(n)}/n)^k)
      \ =\ \pr(\tau_1>t,\ldots,\tau_k>t)
      \ =\ \me(S_t^k),\qquad k\in\nz.
   \end{equation}
   Since $0\le E_t^{(n)}/n\le 1$ and $0\le S_t\le 1$ the convergence (\ref{momconv}) of moments
   implies the convergence $E_t^{(n)}/n\to S_t$ in distribution as
   $n\to\infty$. In order to show that $N_t^{(n)}/n\to S_t$ in
   distribution as $n\to\infty$ it remains to verify that $I_t^{(n)}/n\to 0$
   in distribution as $n\to\infty$. In the following it is verified that
   the latter convergence
   even holds in $L^1$. Each internal branch is generated by a collision.
   Thus, if $C_n$ denotes the total number of collisions, the inequality
   $\me(I_t^{(n)})\le\me(C_n)$ holds. The assumption that the coalescent
   has dust ensures that $C_n/n\to 0$ in $L^1$ by Lemma 4.1 of
   \cite{freundmoehle}. Thus, $I_t^{(n)}/n\to 0$ in $L^1$ as $n\to\infty$,
   which yields the desired convergence $N_t^{(n)}/n\to S_t$ in distribution
   as $n\to\infty$. Thus, the convergence of the one-dimensional distributions
   is established.

   Let us now turn to the proof of the convergence in $D_{[0,1]}[0,\infty)$.
   Let $(T_t^{(n)})_{t\ge 0}$ and $(T_t)_{t\ge 0}$ denote the
   semigroups of $(N_t^{(n)}/n)_{t\ge 0}$ and $(S_t)_{t\ge 0}$
   respectively. By Ethier and Kurtz \cite[p.~172, Theorem 2.11]{ethierkurtz}, applied
   with state spaces $E:=[0,1]$ and $E_n:=\{j/n\,:\,j\in\{1,\ldots,n\}\}$,
   $n\in\nz$, and with the maps $\eta_n:E_n\to E$ and
   $\pi_n:B(E)\to B(E_n)$ defined via $\eta_n(x):=x$ for all $x\in E_n$
   and $\pi_nf(x):=f(x)$ for all $f\in B(E)$ and all $x\in E_n$,
   it suffices to verify that for every $t\ge 0$ and $f\in C(E)$,
   $$
   \lim_{n\to\infty}\sup_{x\in E_n}|T_t^{(n)}\pi_nf(x)-\pi_nT_tf(x)|\ =\ 0.
   $$
   For $f\in C(E)$ and $x\in E_n$ we have
   $T_t^{(n)}\pi_nf(x)=\me(\pi_nf(N_{s+t}^{(n)}/n)\,|\,N_s^{(n)}/n=x)
   =\me(f(N_t^{(nx)}/n))$ and $\pi_nT_tf(x)=T_tf(x)=\me(f(xS_t))$. Thus,
   we have to verify that
   $\lim_{n\to\infty}\sup_{x\in E_n}|\me(f(N_t^{(nx)})/n)-\me(f(xS_t))|=0$.
   Since the polynomials are dense in $C(E)$ it suffices to verify the latter
   equation for monomials $f(x)=x^k$, so we have to prove that
   $$
   \lim_{n\to\infty}\sup_{x\in E_n}|\me((N_t^{(nx)})^k)/n^k - x^k\me(S_t^k)|\ =\ 0,
   \qquad k\in\nz, t\ge 0.
   $$
   Using the decomposition $N_t^{(nx)}=E_t^{(nx)}+I_t^{(nx)}$ and
   the facts that $I_t^{(nx)}\le C_{nx}\le C_n$ and that $C_n/n\to 0$
   in $L^1$ (see Lemma 4.1 of \cite{freundmoehle}), it suffices to show that
   $$
   \lim_{n\to\infty}\sup_{x\in E_n}
   |\me((E_t^{(nx)})^k)/n^k-x^k\me(S_t^k)|\ =\ 0.
   $$
   By (\ref{moments}), $\me((E_t^{(nx)})^k)
   =\sum_{j=1}^k (nx)_jS(k,j)\pr(\tau_{nx,1}>t,\ldots,\tau_{nx,j}>t)$.
   Thus, it suffices to verify that
   $$
   \lim_{n\to\infty}
   \sup_{x\in E_n}\bigg|
   \frac{(nx)_k}{n^k}\pr(\tau_{nx,1}>t,\ldots,\tau_{nx,k}>t)-x^k\me(S_t^k)
   \bigg|\ =\ 0.
   $$
   Since $(nx)_k/n^k\to x^k$ as $n\to\infty$ uniformly on $[0,1]$ it remains
   to prove that
   $$
   \lim_{n\to\infty}\sup_{x\in E_n}
   x^k|\pr(\tau_{nx,1}>t,\ldots,\tau_{nx,k}>t)-\me(S_t^k)|\ =\ 0.
   $$
   This is seen as follows. Choose a sequence
   $(\varepsilon_n)_{n\in\nz}$ satisfying $\varepsilon_n\to 0$ and
   $n\varepsilon_n\to\infty$ (for example $\varepsilon_n:=n^{-1/2}$)
   and distinguish the two cases
   $x\in E_n\cap [0,\varepsilon_n]$ and $x\in E_n\cap (\varepsilon_n,1]$.
   Clearly,
   $$
   \sup_{x\in E_n\cap [0,\varepsilon_n]} x^k
   |\pr(\tau_{nx,1}>t,\ldots,\tau_{nx,k}>t)-\me(S_t^k)|
   \ \le\ 2\varepsilon_n^k\ \to\ 0,\qquad n\to\infty.
   $$
   Moreover, since $p_{t,k}(m):=\pr(\tau_{m,1}>t,\ldots,\tau_{m,k}>t)$ is
   non-increasing in $m$ ($\ge k$), it follows for all $n\in\nz$ with
   $n\varepsilon_n\ge k$ that
   \begin{eqnarray*}
      \sup_{x\in E_n\cap (\varepsilon_n,1]}
      x^k|\pr(\tau_{nx,1}>t,\ldots,\tau_{nx,k}>t)-\me(S_t^k)|
      & \le & \sup_{x\in E_n\cap (\varepsilon_n,1]}
      |p_{t,k}(nx)-\me(S_t^k)|\\
      & \le & p_{t,k}(\lfloor n\varepsilon_n\rfloor) - \me(S_t^k)
      \ \to\ 0
   \end{eqnarray*}
   as $n\to\infty$. The proof of part a) is complete.

   \vspace{2mm}

   {\bf Proof of part b).} We have to verify that
   $(L_t^{(n)}/n)_{t\ge 0}$ converges in $D_{[1,\infty]}[0,\infty)$
   to $(1/S_t)_{t\ge 0}$ as $n\to\infty$. Define
   $\varphi:D_{[0,1]}[0,\infty)\to D_{[1,\infty]}[0,\infty)$ via
   $\varphi(x):=(1/x_t)_{t\ge 0}$ for all $x=(x_t)_{t\ge 0}\in
   D_{[0,1]}[0,\infty)$ with the convention $1/0:=\infty$.
   Since the transformation $\varphi$ is continuous
   we will (equivalently) verify that $(n/L_t^{(n)})_{t\ge 0}$ converges
   in $D_{[0,1]}[0,\infty)$ to $(S_t)_{t\ge 0}$ as $n\to\infty$.

   Let $F_n:=\{n/j\,:\,
   j\in\{n,n+1,\ldots\}\}\cup\{0\}$ and $F:=[0,1]$ denote the state
   spaces and $(U_t^{(n)})_{t\ge 0}$ and $(U_t)_{t\ge 0}$
   the semigroups of $(n/L_t^{(n)})_{t\ge 0}$ and $(S_t)_{t\ge 0}$
   respectively. Define $\pi_n:B(F)\to B(F_n)$ via
   $\pi_n f(x):=f(x)$, $f\in B(F)$, $x\in F_n$.
   By Ethier and Kurtz \cite[p.~172, Theorem 2.11]{ethierkurtz} it suffices
   to verify that for all $t\ge 0$ and all $f\in C(F)$,
   $$
   \lim_{n\to\infty}\sup_{x\in F_n}|U_t^{(n)}\pi_nf(x)-\pi_nU_tf(x)|
   \ =\ 0.
   $$
   For $f\in C(F)$ and $x\in F_n$ we have
   $U_t^{(n)}\pi_nf(x)=\me(\pi_n f(n/L_{s+t}^{(n)})\,|\,n/L_s^{(n)}=x)
   =\me(f(n/L_t^{(n/x)}))$ and $\pi_n U_tf(x)=U_tf(x)=\me(f(xS_t))$. Thus, we
   have to verify that
   $$
   \lim_{n\to\infty}\sup_{x\in F_n}|\me(f(n/L_t^{(n/x)}))-\me(f(xS_t))|
   \ =\ 0.
   $$
   Since the polynomials are dense in $C(F)$ it suffices to verify
   the latter equation for monomials $f(x)=x^k$, so we have to prove that
   $\lim_{n\to\infty}\sup_{x\in F_n}|\me((n/L_t^{(n/x)})^k)-x^k\me(S_t^k)|
   =0$ for all $t\ge 0$ and $k\in\nz$. In the following it is even shown that
   \begin{equation} \label{even}
      \lim_{n\to\infty}\sup_{x\in [0,1]}
      |\me((n/L_t^{(\lfloor n/x\rfloor)})^k)-x^k\me(S_t^k)|
      \ =\ 0,\qquad t\ge 0, k\in\nz,
   \end{equation}
   where $\lfloor n/x\rfloor:=\max\{z\in\gz\,:\,z\le n/x\}$.

   For $m\in\nz$, $t\ge 0$ and $y\in (0,1]$, it follows by
   duality (Theorem \ref{duality} applied with $i:=\lceil m/y\rceil:=
   \min\{z\in\gz\,:\,z\ge m/y\}$ and $j:=m$)
   \begin{eqnarray*}
      \pr(m/L_t^{(m)}\le y)
      & = & \pr(L_t^{(m)}\ge m/y)
      \ = \ \pr(L_t^{(m)}\ge \lceil m/y\rceil)\\
      & = & \pr(N_t^{(\lceil m/y\rceil)}\le m)
      \ = \ \pr\bigg(\frac{N_t^{(\lceil m/y\rceil)}}{\lceil m/y\rceil} \le \frac{m}{\lceil m/y\rceil}\bigg).
   \end{eqnarray*}
   Since $N_t^{(m)}/m\to S_t$ in distribution as $m\to\infty$ by part a)
   of Theorem \ref{theo1}, which is already proven, and since
   $\lim_{m\to\infty} m/\lceil m/y\rceil = y$, we conclude that
   $\lim_{m\to\infty}\pr(m/L_t^{(m)}\le y)=\pr(S_t\le y)$, if $y\in (0,1]$
   is a continuity point of the distribution function of $S_t$. The point
   $y=0$ has to be treated separately. For all $m\in\nz$ we have
   $\pr(m/L_t^{(m)}\le 0)=\pr(L_t^{(m)}=\infty)=
   \lim_{n\to\infty}\pr(L_t^{(m)}\ge n)=\lim_{n\to\infty}\pr(N_t^{(n)}\le m)
   =\lim_{n\to\infty}\pr(N_t^{(n)}/n\le m/n)=\pr(S_t\le 0)$, if $y=0$ is a
   continuity point of the distribution function of $S_t$. The pointwise
   convergence of the distribution functions implies the convergence
   $m/L_t^{(m)}\to S_t$ in distribution as $m\to\infty$.

   Fix $x\in (0,1]$. Replacing $m$ by $\lfloor n/x\rfloor$ it follows by
   an application of Slutzky's theorem that $n/L_t^{(\lfloor n/x\rfloor)}
   \to xS_t$ in distribution as $n\to\infty$. This convergence obviously
   holds as well for $x=0$ with the convention $n/0:=\infty$.
   Noting that the map $y\mapsto y^k$ is bounded and continuous on $[0,1]$
   we conclude that
   $$
   \lim_{n\to\infty}\me((n/L_t^{(\lfloor n/x\rfloor)})^k)
   \ =\ x^k\me(S_t^k),\qquad t\ge 0, k\in\nz, x\in [0,1].
   $$
   In order to see that this pointwise convergence holds even
   uniformly for all $x\in [0,1]$ we proceed as follows.
   Fix $t\ge 0$, $k\in\nz$ and $n\in\nz$. By the pathwise construction
   of the fixation line, we have $L_t^{(1)}\le L_t^{(2)}\le\cdots$. It
   follows that the map $x\mapsto \me((n/L_t^{(\lfloor n/x\rfloor)})^k)$ is
   non-decreasing on $[0,1]$. Clearly, the limiting map
   $x\mapsto x^k\me(S_t^k)$ is non-decreasing,
   bounded and continuous on $[0,1]$. Thus, the pointwise convergence
   holds even uniformly for all $x\in [0,1]$. Note that the proof that
   this pointwise convergence holds even uniformly works
   the same as the proof that pointwise convergence of distribution functions
   is uniform if the limiting distribution function is continuous
   (P\'olya \cite[Satz I]{polya}). Therefore, (\ref{even}) is
   established. The proof is complete.\hfill$\Box$
\end{proof}

\subsection{Appendix} \label{appendix}
\setcounter{theorem}{0}
We establish a sort of duality relation for
generalized Stirling numbers. Let $a,b,r\in\rz$ and suppose
that $t\notin\{0,a,2a,3a,\ldots\}$ such that we can define
$$
q_{ij}\ :=\ \frac{(t-r|b)_j}{(t|a)_i}S(i,j),\qquad i,j\in\nz_0,
$$
where $(t|a)_i:=\prod_{k=0}^{i-1}(t-ak)$ and the coefficients
$S(i,j):=S(i,j;a,b,r)$ are the generalized Stirling
numbers as defined in \cite{hsushiue}. The recursion
$S(i+1,j)=S(i,j-1)+(jb-ia+r)S(i,j)$ for the generalized Stirling numbers
(see \cite[Theorem 1]{hsushiue}) obviously transforms into the recursion
\begin{equation} \label{qrec}
   q_{i+1,j}\ =\ \frac{t-r-(j-1)b}{t-ia}q_{i,j-1} +
   \frac{jb-ia+r}{t-ia}q_{ij}
\end{equation}
for the quantities $q_{ij}$ ($q_{i,-1}:=0$).
Note that $\sum_{j=0}^\infty q_{ij}=1$ for all $i\in\nz_0$.
For $i,j\in\nz_0$ we define $q_{i,\le j}:=\sum_{k=0}^jq_{ik}$.
\begin{lemma} \label{applemma}
   Fix $j\in\nz_0$ and suppose that the limit $\lim_{k\to\infty}q_{k,\le j}$
   exists. Then, for all $i\in\nz_0$,
   \begin{equation} \label{dualrelation1}
      q_{i,\le j} - \lim_{k\to\infty} q_{k,\le j}
      \ =\ \sum_{k=i}^\infty \frac{t-r-jb}{t-ka} q_{kj}.
   \end{equation}
\end{lemma}
\begin{remark}
   At a first glance Lemma \ref{applemma} looks somewhat technical and does
   not seem to have many applications. Indeed, in general it seems to be not straightforward to verify the existence of
   the limit $\lim_{k\to\infty}q_{k,\le j}$ and to determine this limit
   (if it exists). However, for particular parameter
   choices (for instance for $a\le 0$, $b>0$, $r=0$ and $t>0$ an integer
   multiple of $b$), the $q_{ij}$ turn out to be non-negative. In this case there
   exists a random variable $K_i$ with distribution $\pr(K_i=j)=q_{ij}$,
   $j\in\nz_0$.
   Based on the recursion (\ref{qrec}) the sequence $K:=(K_i)_{i\in\nz_0}$
   can be even constructed such that $K$ is a Markov chain with initial state
   $K_0=0$ satisfying $K_{i+1}-K_i\in\{0,1\}$ for all $i\in\nz_0$.
   Hence, $q_{i,\le j}=\pr(K_i\le j)$ is non-increasing in $i$, which ensures
   the existence of the limit $\lim_{i\to\infty}q_{i,\le j}$. If all states
   $0,1,\ldots,j$ of the chain $K$ are transient, then
   $\lim_{i\to\infty}q_{i,\le j}=0$ and (\ref{dualrelation1}) reduces to
   \begin{equation} \label{dualrelation2}
      \sum_{k=0}^j q_{ik}\ =\ \sum_{k=i}^\infty \frac{t-r-jb}{t-ka}q_{kj}.
   \end{equation}
   Roughly speaking, (\ref{dualrelation2}) is a sort of analytic reformulation
   of a particular Siegmund duality. For typical examples we refer the
   reader to the equality of (\ref{diri1}) and (\ref{diri2}) for the Dirichlet
   coalescent and to the equality of (\ref{pd1}) and (\ref{pd2}) for the
   Poisson--Dirichlet coalescent.
\end{remark}
\begin{proof} (of Lemma \ref{applemma})
   The proof of Lemma \ref{applemma} is purely analytic and
   straightforward. For all $k,i\in\nz_0$ we have
   \begin{eqnarray*}
      q_{ki} - q_{k+1,i}
      & = & q_{ki} - \frac{t-r-(i-1)b}{t-ka}q_{k,i-1} -
   \frac{ib-ka+r}{t-ka}q_{ki}\\
      & = & \frac{t-r-ib}{t-ka}q_{ki} - \frac{t-r-(i-1)b}{t-ka}q_{k,i-1}.
   \end{eqnarray*}
   Summation over all $i\in\{0,\ldots,j\}$ yields
   $$
   q_{k,\le j} - q_{k+1,\le j}
   \ =\ \frac{t-r-jb}{t-ka}q_{kj}.
   $$
   Another summation over all $k\ge i$ yields
   the result.\hfill$\Box$
\end{proof}
\begin{acknowledgement}
   The authors thank Jonas Kukla for helpful discussions and comments.
\end{acknowledgement}

\end{document}